\documentclass[final,leqno]{siamltex}
\newtheorem{remark}{Remark}
\newtheorem{example}{Example}
\usepackage{mathtools}
\usepackage{mathpazo}
\usepackage[mathpazo]{flexisym}
\usepackage{breqn}
\usepackage{graphicx}
\usepackage{caption}
\usepackage{subcaption}
\usepackage{aligned-overset}
\usepackage{tikz}
\usepackage{epstopdf}
\usetikzlibrary{calc}
\usepackage[unicode]{hyperref}
\usepackage{cleveref}

\newcommand{\phat}{\widehat{p}}
\DeclareMathOperator{\Beta}{B}

\newcommand{\pFq}[5]{{}_{#1}{F}_{#2}\left(\begin{matrix}#3\\#4\end{matrix};#5\right)}
\newcommand{\gfrac}[2]{\frac{\Gamma(#1)}{\Gamma(#2)}}
\newcommand{\R}{\mathbb{R}}
\newcommand{\N}{\mathbb{N}}

\newcommand{\V}{\mathbb{V}}
\newcommand{\refel}{\hat{\triangle}}
\newcommand{\T}{\triangle}
\newcommand{\triang}{\mathcal{T}}
\newcommand{\vek}[1]{\underline{#1}}  
\newcommand{\dx}{\;\mathrm{d}x}
\newcommand{\dy}{\;\mathrm{d}y}
\newcommand{\dz}{\;\mathrm{d}z}
\newcommand{\scalarp}[2]{\langle #1,#2 \rangle}

\begin{document}
\title{Recursion formulas for integrated products of Jacobi polynomials}
\author{Sven Beuchler
\thanks{Leibniz University Hannover, Institute of Applied Mathematics (IfAM), Welfengarten 1,
30167 Hannover, Germany and Cluster of Excellence PhoenixD (Photonics, Optics, and Engineering - Innovation Across
Disciplines), Leibniz University Hannover, Germany, (\href{mailto:beuchler@ifam.uni-hannover.de}{beuchler@ifam.uni-hannover.de})}
\and Tim Haubold
\thanks{
Leibniz University Hannover, Institute of Applied Mathematics (IfAM), Welfengarten 1,
30167 Hannover, Germany, (\href{mailto:haubold@ifam.uni-hannover.de}{haubold@ifam.uni-hannover.de})}
\and Veronika Pillwein
\thanks{
Johannes-Kepler-University Linz, Research Institute for Symbolic Computation (RISC), Altenberger Stra{\ss}e 69, 4040 Linz, Austria,
(\href{mailto:veronika.pillwein@risc.uni-linz.at}{veronika.pillwein@risc.uni-linz.at})}} 
\date{\today}
\maketitle

\begin{abstract}
From the literature it is known that orthogonal polynomials as the Jacobi polynomials can be expressed
by hypergeometric series.
In this paper, the authors derive several contiguous relations for terminating multivariate hypergeometric
series. With these contiguous relations one can prove several recursion formulas of those series. 
This theoretical result allows to compute integrals over products of Jacobi
polynomials in a very efficient recursive way.
Moreover, the authors present an application to numerical analysis where it can be used
in algorithms which compute the approximate solution of boundary value problem of partial differential equations
by means of the finite elements method (FEM).
With the aid of the contiguous relations, the approximate solution can be computed much faster than using numerical integration.
A numerical example illustrates this effect.
\end{abstract}



\section{Introduction}


Finite element methods (FEM) are popular and versatile methods to compute
approximate solutions of partial differential equations (PDE) on complicated
domains, for which no exact solution is known. The solution is expanded in a
basis that is constructed on a mesh of basic geometric objects, in our case
simplices. The coefficients of these basis functions are computed as solution to
a linear system $K\vek{u}=\vek{f}$ whose entries are the integral of products of
(derivatives of) these basis functions, see
~\cite{Ciarlet,Szabo,QuatValli,Braess}).  If the solution of the underlying PDE
is smooth, the local polynomial degree $p$ is increased in order to improve the
accuracy of the computed FEM-solution in comparison to the exact one,
\cite{Schwab,Dem08}.  However, the computation of the FEM-solution requires a
suitable basis in order to keep the computational cost as low as possible.
In~\cite{beuchler}, a set of basis functions based for the Poisson problem in 2D
was presented, that yields a sparse linear system.  These basis function are
based on classical orthogonal polynomials on a triangle, see
e.g. \cite{Proriol},\cite{Koornwinder} and also \cite{Dubiner}.  The sparsity is
obtained by the orthogonality relations.  This work was later generalized to 3D
as well as other partial differential equations, see~\cite{Zaglmayer} for an
overview on these results.

In all these cases, the sparsity for the respective basis functions was proven
by explicitly evaluating the integrals using difference-differential relations
satisfied by Jacobi polynomials. In the case of tetrahedral elements, these
computations became so involved that they could not be carried out by hand, but
were proven using symbolic computation. The basis for the application of
computer algebra is that Jacobi polynomials $P_n^{(\alpha,\beta)}(x)$ are
holonomic~\cite{Zeil90a,Koutschan09}, i.e., they satisfy certain types of
difference-differential equations in all parameters. For various
representations, it is possible to automatically derive new identities involving
Jacobi polynomials automatically~\cite{kauers13,RISC2949}. These procedures are
rigorous and many come with easy-to-verify certificates.

The outcome of the algorithm are the coefficients $K_{ik}$ of the linear system
(matrix) $K$ computed explicitly as rational functions in the parameters. For a
rational function, it is fairly easy to derive recurrence relations. In order to
find a minimal order recurrence, we have used a Mathematica package by
Kauers~\cite{Kauers} for guessing recurrences.  These are easily proven by
plugging in the actual entries. This guess-and-prove approach is very common
also with more complicated input than rational functions.  There exist
implementations of algorithms for holonomic functions in several computer
algebra systems such as, e.g., \texttt{mgfun} in Maple~\cite{Chyzak98b} or
\texttt{ore}\textunderscore algebra in Sage~\cite{oresage}.


The results in this article show how recurrences, for the integrals of two
multivariate orthogonal polynomials and $H^1$ basis functions on a triangle, can
be found directly from a multivariate hypergeometric series representation
$_pF_q $ of the integrals using contiguous relations and terminating $_pF_q$
identities. This kind of techniques is well known since the works of Gau\ss~and
Euler, see e.g. \cite{Askey} or \cite{Rainville}, but are seldom applied to
multivariate series.  Usually one would refer to the works of Wilson
\cite{Wilson78} or would use symbolic computation to derive such recursion
formulas, but a direct derivation approach happened to be more insightful, not
only for the proof, but for the general structure of such series as well.
Multivariate hypergeometric series are more difficult to handle than general
$_pF_q$ series. In some cases one can reduce a multivariate series to a general
one, by using well known summation theorems like the Pfaff-Saalsch\"utz,
Dougall's Summation or Whipple's transformation, see e.g. \cite{Bailey35} or
\cite{Askey}, but usually there is no transformation, which holds for all
parameter configurations, which we are interested in.  A broad overview of
convergence theory for multivariate series can be found e.g. in \cite{Exton},
but does not need to be applied here, since our series are based on Jacobi
polynomials, which are terminating series in itself.  The notation which will be
used goes back to Burchnall and Chaundy \cite{Burchnall}.

For the first time, this sheds light on the underlying
structure of the integral values and furthermore the sparsity pattern of finite element matrices can be read out from the coefficients of the recursion formulas.
To our knowledge, these identities are unknown and interesting in their own right.
For the community in numerical analysis theorem \ref{MainTheo} is
the most important result stating how the nonzero entries of $K$ can be computed recursively in optimal
arithmetical complexity. It is a consequence of the main result of this article, theorem~\ref{Main}. 
The results can be extended to an arbitrary simplex or basis functions in the function spaces $H(\operatorname{curl})$ and $H(\operatorname{div})$.
Optimal complexity was first achieved for finite element matrices not so long ago by \cite{Ainsworth1} for basis functions based on Bernstein polynomials, but the resulting element
matrices were dense. One can transform the basis functions based on Jacobi polynomials to a basis based on Bernstein polynomials, see e.g. \cite{Ainsworth2},\cite{Ainsworth3}, 
but this transformation loses optimal arithmetical complexity for the assembly, though it has other useful properties, which will not be discussed here.
For element mass matrices, based on Bernstein polynomials, one can use a block recursive structure as can be seen in \cite{Kirby}, 
which results in an efficient inversion technique. Furthermore recursion formulas for the orthogonal polynomials on a triangle have already been computed as can be seen in e.g. in \cite{Xu} and \cite{Townsend1}.
But to the knowledge of the authors, there are only a few
publications in which the matrix entries $K_{ik}$ of a FEM-Matrix are computed
completely recursively, see \cite{PillweinPamm} for a special case. Part of the here shown recursion formulas were first
published in \cite{Pamm}. This paper generalizes the recursion formulas and presents a proof for those relations.

This paper is organized as follows:
In section \ref{ch2}, the authors introduce hypergeometric series and Jacobi polynomials and give an overview on several known identities and notations that are used throughout the paper.
A short motivation of the background from FEM is given in section~\ref{sec:fem}.
The main theoretical results are formulated in section \ref{ch4}.
Finally, some algorithmic aspects and numerical experiments are presented in section \ref{ch5}.

Throughout this paper, the indices $n,m,i,j,k,l$ denote natural numbers. $P^{(\alpha,\beta)}_n$ is the $n$-th Jacobi Polynomial with indices $\alpha,\beta.$ The set $\rho,\delta$ usually denotes another set of Jacobi polynomial indices.
$F$ stands for some kind of (generalized or multivariate) hypergeometric series and $I$ for the exact value of an integral.
%
%
%
 \section{\label{ch2}Introduction to special functions}
 \subsection{Hypergeometric series}
 For $a \in \R$ let
 \begin{equation*}
  (a)_n = (a)(a+1) \dots (a+n-1) = \gfrac{a+n}{a}
 \end{equation*}
be the Pochhammer symbol or rising factorial and $\Gamma()$ denote the Gamma function~\cite{NIST:DLMF}, which is given by $\Gamma(n) = (n-1)!$ for $n \in \mathbb{N}$.
 \begin{definition}[Gaussian hypergeometric series]
  For $a,b,c$ arbitrarily the series
  \begin{equation}\label{Gauss}
   \pFq{2}{1}{a,b}{c}{x} = \sum_{n=0}^\infty \frac{(a)_n (b)_n}{(c)_n} \frac{x^n}{n!}
  \end{equation}
is called (Gaussian) hypergeometric series.\\
 \end{definition}
 The Gaussian hypergeometric series converges absolutely for $\operatorname{Re}(c-a-b) > 0$, see \cite{Askey} or \cite{Rainville}.
 Classical orthogonal polynomials can be expressed as hypergeometric series, see e.g.~\cite{Rainville}. For Jacobi polynomials
 in particular, we have
\begin{equation}\label{JacPol}
 P_n^{(\alpha,\beta)}(x) = \frac{(\alpha+1)_n}{n!} \pFq{2}{1}{-n,n+\alpha+\beta+1}{\alpha+1}{\frac{(1-x)}{2}}.
\end{equation}
For $x=1$ series \eqref{Gauss} can be summed and written in closed form.
\begin{theorem}[Gau\ss$\lbrack1812\rbrack$]
For $\operatorname{Re}(c-a-b) >0$
\begin{equation}\label{GaussSum}
 \pFq{2}{1}{a,b}{c}{1} =\sum_{n=0}^\infty \frac{(a)_n (b)_n}{(c)_n n!} = \frac{\Gamma(c)\Gamma(c-a-b)}{\Gamma(c-a)\Gamma(c-b)}. 
\end{equation}
 \end{theorem}
This can be proven by using an integral representation due to Euler, see \cite{Askey}.
If $a$ (or $b$) is a negative integer, the identity simplifies even more:
\begin{corollary}[Chu-Vandermonde]
Let $a = -m$ with $m\in\N$. Then
\begin{equation*}
 \sum_{n=0}^\infty \frac{(-m)_n (b)_n}{(c)_n n!} = \frac{(c-a)_m}{(c)_m}.
\end{equation*}
 \end{corollary}
A generalized version of the ${}_2F_1()$ is called a generalized hypergeometric series. 
\begin{definition}
 Let $(a_i),(b_j), i = 1,\dots p, j = 1,\dots q$. Then the series
 \begin{equation}\label{genHyp}
  \pFq{p}{q}{a_1,a_2, \dots ,a_p}{b_1,b_2,\dots b_q}{x} = \sum_{n=0}^\infty \frac{(a_1)_n(a_2)_n \dots (a_p)_n}{(b_1)_q (b_2)_q \dots (b_n)_q } \frac{x^n}{n!}
 \end{equation}
is called a generalized hypergeometric series. \end{definition}\\
The series $\pFq{p}{q}{a_1,a_2, \dots ,a_p}{b_1,b_2,\dots b_q}{x}$ converges absolutely for all $x$ if $p<q$ or for $|x|<1$ if $p=q+1$, see e.g. \cite{Askey}.
Higher classes of polynomials, e.g. polynomials of the \textit{Hahn} class (see \cite{Askey}), can be described by those series.\\
For the special case $\pFq{3}{2}{-m,a,b}{c,1+a+b-c-m}{1}$, $m\in\N$, a similar summation to \eqref{GaussSum} can be found.
\begin{theorem}[Pfaff-Saalsch\"utz]
Let $m \in \N$, then  
 \begin{equation}\label{Pfaff}
 \sum_{n=0}^\infty \frac{(-m)_n (a)_n (b)_n}{(c)_n (1+a+b-c-n)_n n!} = \frac{(c)_m (c-a-b)_m}{(c-a)_m (c-b)_m}. 
\end{equation}
\end{theorem}
This can be proven by equating the coefficients of a transformation for the ${}_2F_1$ which again is due to Euler, see \cite{Askey} or \cite{Rainville}. Important for summability is the fact that the series is balanced or Saalsch\"utzian.
\begin{definition}
 A generalized hypergeometric series \eqref{genHyp} is called $s$-balanced if $x=1$ and \begin{equation*}a_1 + a_2 + \dots a_p = b_1 + b_2 + \dots b_q + s.\end{equation*} \\The case $s=1$ is also called Saalsch\"utzian. \\
\end{definition}
There are some summation theorems like Chu-Vandermonde or the Pfaff-Saalsch\"utz theorem for $p,q$ greater than $3,2$, but they are usually more restrictive. For example for a balanced ${}_7F_6$ series, there holds Dougall's summation, but the series needs to fulfil additional properties, see e.g. \cite{Askey},\cite{Bailey35} for more information. 
\subsection{Contiguous Relations}
Recurrence relations can be proven by using the contiguous relations of the hypergeometric series. We will briefly summarize the basic ideas as can be found in the book of Rainville \cite{Rainville}. Alternative and/or equivalent methods for more general hypergeometric series, can be found for example
in the book of Andrews et al. \cite{Askey} or in the work of Bailey \cite{Bailey54} and Wilson \cite{Wilson78}.\\
Let $F \coloneqq \pFq{2}{1}{a,b}{c}{x}$ be a Gaussian hypergeometric series. Taking the derivative yields
\begin{align}\label{GaussDeriv}
 \frac{\partial}{\partial x} F &= \sum_{n=0}^\infty \frac{(a)_n (b)_n}{(c)_n} \frac{n x^{n-1}}{n!} =  \sum_{n=1}^\infty \frac{(a)_n (b)_n}{(c)_n} \frac{x^{n-1}}{(n-1)!}
=  \sum_{n=0}^\infty \frac{(a)_{n+1} (b)_{n+1}}{(c)_{n+1}} \frac{x^{n}}{n!} = \frac{a b}{c} \pFq{2}{1}{a+1,b+1}{c+1}{x}.
\end{align}
For the ease of notation when stating the contiguous relations, we follow common practice~\cite{Askey} and write
\begin{align*}
 F(a+) &\coloneqq \pFq{2}{1}{a+1,b}{c}{x},  & F(a-) &\coloneqq \pFq{2}{1}{a-1,b}{c}{x},\\
  F(b+) &\coloneqq \pFq{2}{1}{a,b+1}{c}{x},  & F(b-) &\coloneqq \pFq{2}{1}{a,b-1}{c}{x},\\ 
  F(c+) &\coloneqq \pFq{2}{1}{a,b}{c+1}{x},  & F(c-) &\coloneqq \pFq{2}{1}{a,b}{c-1}{x}.
\end{align*}
Therefore \eqref{GaussDeriv} can be written as $\frac{\partial}{\partial x} F =  \frac{a b}{c} F(a+,b+,c+).$ Now define the differential operator $\theta_x = x \frac{\partial}{\partial x}$, also known as Euler operator. It is applied as follows
\begin{equation}
\label{Gaussa}
 (\theta_x +a) F = \sum_{n=0}^\infty \frac{(a)_n (b)_n}{(c)_n} \frac{(a+n) x^{n}}{n!}= \sum_{n=0}^\infty \frac{(a)_{n+1} (b)_n}{(c)_n} \frac{x^{n}}{n!}
= a \sum_{n=0}^\infty \frac{(a+1)_n (b)_n}{(c)_n} \frac{x^{n}}{n!}= a F\left(a+\right).
\end{equation}            
Analogously the formulas
\begin{align}
 \label{Gaussb} (\theta_x +b) F &= b F(b+) \quad \text{and}\\
 (\theta_x +c-1) F &= (c-1) F(c-) \quad \text {can be proven.}
\end{align}
From \eqref{Gaussa} and \eqref{Gaussb} follows
\begin{equation*}
 (a-b) F = a F(a+) - b F(b+),
\end{equation*}
this is one of the 15 contiguous relations\footnote{Although this reduces to 9 relations, if one takes symmetry into account, see \cite{Askey}} for the ${}_2F_1$ and its six contiguous functions. The other can be obtained by the same kind of straight forward computations. For a complete list, see e.g. \cite{Askey}, \cite{Bailey35}, \cite{Rainville} or \cite{Slater}.\\ Many of the important recurrence relation between Jacobi polynomials can be proven by using the contiguous relations of the ${}_2 F_1$. 
\subsection{Jacobi Polynomials}
The Jacobi polynomials \eqref{JacPol} are the polynomials, which are orthogonal on $[-1,1]$ with respect to the weight function $w(x) \coloneqq (1-x)^\alpha (1+x)^\beta, \alpha,\beta>-1$. They can either be given by
\begin{equation*}
 P_n^{(\alpha,\beta)}(x) = \frac{1}{(x-1)^\alpha (x+1)^\beta 2^n n!} \frac{\partial}{\partial x}^n (x-1)^{n+\alpha}(x+1)^{n+\beta} , 
\end{equation*}
which is the Rodrigues formula or by the hypergeometric representation \eqref{JacPol}.
Furthermore the property
\begin{equation}\label{atob}
 P_n^{(\alpha,\beta)}(x) = (-1)^n P_n^{(\beta,\alpha)}(x)
\end{equation}
yields the representation
\begin{equation}\label{JacPol2}
 P_n^{(\alpha,\beta)}(x) = \frac{(-1)^n (1+\beta)_n}{n!} \pFq{2}{1}{-n, n+\alpha+\beta+1}{\beta+1}{\frac{1+x}{2}}.
\end{equation}
We refer to the standard literature for more information on these properties,
e.g. \cite{Askey}, \cite{Rainville},\cite{Szego}.  Since the Jacobi polynomials
are orthogonal polynomials they satisfy a three term recurrence relation. It is
given by,
\begin{equation}\label{3term}
\begin{aligned}
 2n (\alpha+&\beta+n)(\alpha+\beta+2n-2)P_n^{(\alpha,\beta)}(x) \\ =& \left[(\alpha+\beta+2n-1)(\alpha^2 - \beta^2) + x(\alpha+\beta+2n-2)_3\right] P_{n-1}^{(\alpha,\beta)}(x)\\ &- 2(a+n-1)(\beta+n-1)(\alpha+\beta+2n) P_{n-2}^{(\alpha,\beta)}(x).
\end{aligned}
\end{equation}
Next, let us recall several difference equations satisfied by 
Jacobi polynomials,
(see e.g. \cite[Ch. 16]{Rainville}) that are summarized in the following lemma.
\begin{lemma}\label{le:Rek}
The Jacobi polynomials \eqref{JacPol} satisfy
\begin{align}
 (\alpha + \beta +n) P_n^{(\alpha,\beta)}(x) &= (\beta + n) P_{n}^{(\alpha,\beta-1)}(x) + (\alpha+n) P_n^{(\alpha-1,\beta)}(x)\\
 \frac{1}{2} (2+\alpha + \beta +2n) (x-1) P_n^{(\alpha+1,\beta)}(x) &= (n+1) P_{n+1}^{(\alpha,\beta)}(x) - (1+\alpha+n)P_n^{(\alpha,\beta)}(x) \\
 \frac{1}{2} (2+\alpha + \beta +2n) (x+1) P_n^{(\alpha,\beta+1)}(x) &= (n+1) P_{n+1}^{(\alpha,\beta)}(x) + (1+\beta+n) P_n^{(\alpha,\beta)}(x) \\
 (\alpha + \beta +2n) P_n^{(\alpha,\beta-1)}(x) &= (\alpha+\beta+n)P_n^{(\alpha,\beta)}(x) + (\alpha +n) P_{n-1}^{(\alpha,\beta)}(x) \\
 (\alpha + \beta +2n) P_n^{(\alpha-1,\beta)}(x) &= (\alpha+\beta+n)P_n^{(\alpha,\beta)}(x) - (\beta+n) P_{n-1}^{(\alpha,\beta)}(x)\\
 2 P_n^{(\alpha,\beta)}(x) &= (1+x) P_n^{(\alpha,\beta+1)}(x) + (1-x)P_n^{(\alpha+1,\beta)}(x) \\
 P_{n-1}^{(\alpha,\beta)}(x) &= P_{n}^{(\alpha,\beta-1)}(x) - P_{n}^{(\alpha-1,\beta)}(x)  
\end{align}
\end{lemma}
Further relations between different Jacobi polynomials that can be found in~\cite{Rainville} will be introduced if required.

In a high order finite element context it is often useful to use integrated Jacobi polynomials $\hat{P}_n^{(\alpha,0)}(x)$,, which can be written as Jacobi polynomials, 
\begin{equation}\label{IntJac}
 \int_{-1}^x P_{n-1}^{(\alpha,0)}(s) \,ds \eqqcolon \hat{P}_n^{(\alpha,0)}(x) = \frac{2}{n+\alpha-1}P_n^{(\alpha-1,-1)}(x),
\end{equation}
where the Jacobi polynomials with index $\beta = -1$ and $\alpha > -1$ are defined properly, see e.g. \cite{Szego}, as
\begin{equation}\label{NotIntJac}
 P_n^{(\alpha,-1)}(x) = \frac{1+x}{2} \frac{n+\alpha}{n} P_{n-1}^{(\alpha,1)}(x).
\end{equation}
Integrated Legendre polynomials can be defined similarly,
\begin{equation}\label{NotIntLeg}
 \int_{-1}^{x} P_{n-1}^{(0,0)}(s) \,ds = \hat{P}_n^{(0,0)} (x)= \frac{2}{n-1}P_n^{(-1,-1)}(x) = \frac{x^2-1}{2(n-1)} P_{n-2}^{(1,1)}(x).
\end{equation}
\section{Background from the Finite Element Method (FEM)}\label{sec:fem}
\paragraph{Variational formulation and the function space $H^1$}
In this paper, we investigate the following problem in variational formulation:
For a given bounded domain $\Omega\subset \R^d$ 
find $u$ in the Sobolev space $H^1$ 
such that
\begin{equation}
  \label{blf-def}
  a(u,v):=\int_{\Omega} \mu \nabla u \cdot \nabla v + \int_{\Omega} \kappa \, u\; v=\int_{\Omega} f  v=: F(v)
\quad\forall v \in  H^1
\end{equation}
holds. For ease of notation, we assume Neumann boundary conditions.
The bilinear form $a(\cdot,\cdot)$ and the linear form $F(\cdot)$ are
well-defined and bounded for $f\in[L_2(\Omega)]$ and
$\mu,\kappa\in[L_\infty (\Omega)]$ with $\mu>0$ and $\kappa,\mu$ are
assumed to be piecewise constant.  The variational formulation is
well-defined for square-integrable vector-valued functions $u: \Omega
\rightarrow \R^3$ with square-integrable gradient. We denote the
according function space by \begin{align} H^1 := \{ u \in
  [L_2(\Omega)] \: : \: \nabla u \in [L_2(\Omega)]^d\}, \quad d=2,3. \end{align}

The variational formulation~\eqref{blf-def} is obtained from the discretization of the
reaction-diffusion equation $-\nabla \cdot (\mu \nabla u)+\kappa u=f$ by multiplication with a test function $v$, integration
over the domain $\Omega$ and applying Green's formula to the second order part. We refer the
interested reader to~\cite{Evans,MR0350177,MR1814364} for more informations concerning this topic.

For complicated geometries $\Omega$ and real-life data it is not possible to solve the
equations~\eqref{blf-def} analytically.  The finite element method (FEM) provides
a general method for the numerical solution of partial differential equations.  It is
based on the variational formulation of the underlying PDE and provides a profound
analysis.

\paragraph{Finite element discretization of $H^1$}
Galerkin methods such as, e.g., FEMs are among the most powerful
methods for the solution of boundary value problems of the
form~\eqref{blf-def}.  The Galerkin approximation relies on the
orthogonal projection of the implicitly given solution $u$
in~\eqref{blf-def} onto a $N$-dimensional subspace $\V_{N}
\subset H^1$ with respect to the bilinear form $a(\cdot,\cdot)$.
Therefore, we construct a sequence of finite dimensional spaces $\V_N
\subset H^1$ and consider the solution of
\eqref{blf-def} on $\V_N$ (see e.g.~\cite{Ciarlet,Szabo} or the
textbooks \cite{QuatValli,Braess}), namely
\begin{equation}
  \label{discr-probl}
  \text{Find $u_N \in \V_N$ such that} \qquad a(u_N,v_N)=F(v_N) \quad\forall v_N\in \V_N.
\end{equation}
The finite element method provides a special construction of these
discrete spaces $\V_N$ by piecewise polynomial functions on
an admissible subdivision (see~\cite{Ciarlet}) $\mathcal{T}_h$
of~$\Omega$ into simplices $\tau_s$ with $s=1,...,nel$, i.e.,
\begin{equation}\label{eq2p3}
\V_{N} := \{ u \in H^1 \: : \: u|_{\tau_s} \in {\mathcal{P}^p(\tau_s)} \: \forall \tau_s \in \triang_h \},
\end{equation}
where $\mathcal{P}^p$ is the space of all polynomials defined on $\tau_s$ of maximal total
degree~$p$.  The elements $\tau_s$ are chosen such that $\kappa$ and $\mu$ are constant on
the elements.  In $hp$-finite element methods the polynomial degree $p$ can vary on each
element $\tau_s$ which provides extraordinary fast convergence of the finite element
method with respect to the number of degrees of freedom $N = \dim (\V_N)$, see
e.g.~\cite{Schwab}.  This is crucial for the solution of real world problems of the
form~\eqref{blf-def}.

Since the space $\V_N$ is finite dimensional, the space is equipped with a row vector of
basis functions $[\Psi]:=[\psi_1,\ldots,\psi_N]$.  
The basis functions $\psi_j$ are chosen such that they have local support (see e.g.~\cite{Ciarlet}).

Then using the ansatz $u_N(x)=\sum_{i=1}^N u_i \psi_i(x)$ and setting $v = \psi_j$ for
$j=1,...,N$ in~\eqref{discr-probl} the problem becomes equivalent to solving the following
system of $N$ linear algebraic equations
\begin{equation}
  \label{eq2p5}
 \text{Find a coefficient vector $\vek{u}:=[u_i]_{i=1}^{N} \in \R^N$ s.t.} \quad   K \vek{u}=\vek{f}  
\quad\text{with} \begin{array}{rl}  K=\left[a(\psi_j,\psi_i)\right]_{i,j=1}^N \in \R^{N\times N} & \;\text{system matrix,} \\ 
\vek{f}=\left[F(\psi_j) \right]_{j=1}^N \in \R^{N} &\;\text{right hand side vector.}  \end{array}
\end{equation} 

Note that the matrix $K$ depends on the choice of the basis functions.

\paragraph{Efficient solution of algebraic system}
In practical problems, the dimension $N$ usually becomes very large ($\succeq 10^6$).
Iterative methods as the preconditioned GMRES method or the preconditioned conjugate
gradient-method (pcg-method) for positive definite systems are preferred for the solution
of~\eqref{discr-probl}.  The two main important issues for the fast solution of the system
$K\vek{u}=\vek{f}$ are
\begin{itemize}
\item the fast multiplication $K\vek{u}$,
\item the choice of a good preconditioner $C^{-1}$ for $K$ such that the condition number
  $\kappa(C^{-1}K)$ becomes small,
\end{itemize}
in order to obtain a fast convergence of the iterative solver for the
solution of \eqref{eq2p5}.  If $K$ is a dense matrix, the operation
$K\vek{u}$ requires $N^2$ flops.  If $K$ is a sparse matrix with a
bounded number $c$ of nonzero entries per row, the computational
complexity of the matrix vector-multiplication is bounded by
$cN$. Since $K$ in~\eqref{eq2p5} depends on the choice of the basis
$[\Psi]$ it is essential to choose a basis with as many orthogonality
relations as possible with respect to the bilinear form
$a(\cdot,\cdot)$.  The choice of the basis heavily influences the
properties of the matrix $K$:
\begin{itemize} 
\item the local support of finite element basis functions yields sparse system matrices
  $K$ and hence a cheap matrix vector multiplication $K \vek{u}$
\item the condition number of $K$ and $C^{-1} K$, respectively, stability and less
  iterations in iterative solution methods.
\end{itemize} 
In the lower order version of FEM, i.e., the $h$-version, multigrid solvers are the most
powerful methods for discretizations of boundary value problems of partial differential equations, see \cite{Hack1} and the references therein.  
For $hp$-FEM this strategy is combined with appropriate local smoothers and static condensation.

\paragraph{$hp$-FEM and choice of basis functions} 
In $hp$-FEM, the local polynomial degree $p_s$ on the elements may be large.  Despite of
the local support of the basis functions $[\Psi]$, the local dimension $n_s$ grows as
$\mathcal{O}(p_s^d)$, where $d=2,3$ is the spatial dimension.  
Hence we are interested in a bounded number of nonzero entries in
the system matrix independent of the polynomial degrees.  Let
$[\Phi_s]=[\phi_{i,s}]_{i=1}^{n_s}$ be the set of all basis functions $\psi_j$ with $\supp
\psi_j\cap \tau_s\neq \emptyset$, e.g. $[\Phi_s]=[\Psi] L_s$ with the (boolean) finite
element connectivity matrices~$L_s$.

In finite element methods, the global system matrix $K$ is the result of assembling local
matrices, see~\cite{Ciarlet}.  In our case, one obtains
\[
K=\int_{\Omega} \mu \nabla [\Psi]^\top \cdot \nabla [\Psi]  + \kappa \, [\Psi]^\top  \;  [\Psi] 
=\sum_{s=1}^{nel} \int_{\tau_s} \mu \nabla [\Psi]^\top \cdot \nabla [\Psi]  + \kappa \, [\Psi]^\top  \;  [\Psi].
\]
Together with $[\Phi_s]=[\Psi] L_s$ on $\tau_s$, this implies
\begin{equation}
  \label{eq2p6}
  K=\sum_{s=1}^{nel} L_s^\top K_s L_s
\end{equation}
with the local stiffness and mass matrices
\begin{eqnarray}
  \label{eq2p7}
  K_s&=& \int_{\tau_s}  \left(\mu^{-1} \nabla [\Phi_s]^\top  \cdot \nabla [\Phi_s]
+\kappa  [\Phi_s]^\top \cdot [\Phi_s] \right)  \\
&=& \mu_s \int_{\tau_s} \nabla [\Phi_s]^\top \cdot \nabla [\Phi_s] +
\kappa_s  \int_{\tau_s} [\Phi_s]^\top \; [\Phi_s] =:
\mu_s A_s+\kappa_s M_s \nonumber
\end{eqnarray}
 on the elements $\tau_s$, respectively,
where $\mu_s=\mu\mid_{\tau_s}$ and  $\kappa_s=\kappa\mid_{\tau_s}$ are constants.

Therefore, the sparsity of the matrices $K_s$ in~\eqref{eq2p7} implies sparsity of
 the matrix $K$, cf.~\eqref{eq2p6}.  Our aim is to develop a local polynomial basis $[\Phi_s]$
 such that the matrices $A_s$ and $M_s$ in~\eqref{eq2p7} have a bounded number of nonzero
 entries per row.
 The global basis is the obtained in the usual way, see e.g. \cite{Dem08}.

\paragraph{Model problem} For ease of presentation, we are focusing on the following model
 problem: Let $\T$ denote an arbitrary non degenerated simplex
 $\T\subset\mathbb{R}^2$.  Find a polynomial basis
 $[\Phi]=\left[\phi_i\right]_{i=1}^{n(p)}$ of degree $p$ with $\phi_i: \T\mapsto \mathbb{R}^2$
such that the
 matrices
\begin{eqnarray}
  M&:=&  \left[\int_{\T} \phi_{j} \; \phi_{i} \right]_{i,j=1}^{n}=\int_{\T}[\Phi]^\top \;[\Phi]   \label{eq2p8}\\
  A&:=&  \left[\int_{\T} \nabla \phi_{j} \cdot \nabla \phi_{i}\right]_{i,j=1}^{n}= \int_{\T}\nabla [\Phi]^\top \nabla[\Phi]\nonumber
\end{eqnarray}
have $\mathcal{O}(n)$ nonzero entries.  This basis should be suited for $H^1$ conformity.

\paragraph{Definition of the basis functions}
Let $\T$ denote an arbitrary non-degenerated simplex $\T \subset \R^2$, its set of four
vertices by $\mathcal{V} = \{V_1,V_2,V_3\}$, $V_i \in \R^2$ and $\lambda_1, \lambda_2,
\lambda_3 \in \mathcal{P}^1(\T)$ its barycentric coordinates uniquely defined by $\lambda_{i}(V_j) = \delta_{ij}$.

Using the integrated Jacobi polynomials \eqref{IntJac}, we define the shape functions on the affine triangle
$\T$ with baryzentrical coordinates $\lambda_m(x,y)$,
$m=1,2,3$.
\begin{itemize}
\item The vertex functions are chosen as the usual linear hat functions
 \begin{equation}\nonumber
  \psi_{V,m}(x,y):=\lambda_m(x,y),\quad m=1,2,3.
\end{equation}
Let $\Psi_V^2:=\left[ \psi_{V,1}, \psi_{V,2}, \psi_{V,3}\right]$ be
the basis of the vertex functions.
\item For each edge $E=[e_1,e_2]$,
running from vertex $V_{e_1}$ to $V_{e_2}$, we define
\begin{displaymath}
  \psi_{[e_1,e_2] ,i}(x,y)= \hat{p}_i^{(0,0)}\left( \tfrac{\lambda_{e_2} -
    \lambda_{e_1}}{\lambda_{e_1} +\lambda_{e_2}} \right) \left(
  \lambda_{e_1} + \lambda_{e_2} \right)^i.
\end{displaymath}
By $\Psi_{[e_1,e_2]}=\left[ \psi_{[e_1,e_2],i} \right]_{i=2}^{p}$, we  denote the basis of the edge bubble
functions on the edge $[e_1,e_2]$.
$\Psi_E^2=\left[  \Psi_{[1,2]}, \Psi_{[2,3]}, \Psi_{[3,1]} \right]$ is the
basis of all edge bubble functions.
\item The interior bubbles are defined as
 \begin{equation}
  \label{eq3p6}
  \psi_{ij}(x,y):=g_i(x,y)h_{ij}(x,y),
 \quad i+j\leq p,i\geq 2, j\geq 1,
\end{equation}
where the auxiliary bubble functions $g_i$ and $h_{ij}$ are given by
\[
g_{i}(x,y) := \hat{p}_i^{(0,0)}\left( \tfrac{\lambda_{2} - \lambda_{1}}{\lambda_{1} +\lambda_{2}} \right)\left(\lambda_{1}+\lambda_{2}\right)^i 
\quad\textrm{and}\quad  h_{ij}(x,y):= \hat{p}_j^{(2i-1,0)}(2 \lambda_3 - 1), 
\]
Moreover,
$\Psi^2_I=\left[\psi_{ij} \right]_{i\geq 2,j\geq 1}^{i+j\leq p}$ denotes the basis of all
interior bubbles.
\end{itemize}
Finally, let $\Psi_{\nabla,2}=\left[\Psi^2_V,\Psi^2_E,\Psi^2_I
  \right]$ be the set of all shape functions on $\triangle_s$.

\paragraph{Sparsity Results}
It can be proved, \cite{beuchler}, that the matrices $M$ and $A$ \eqref{eq2p8} have a limited number of nonzeros entries. Usually the stiffness matrices is the more important part, but for the ease of presentation, we will focus on the mass matrix.

In particular, one obtains
\begin{equation}\label{eq3p7}
m_{ij,kl}^{(\T)}=\int_{\T} \psi_{ij}(x,y) \psi_{kl}(x,y)=0 \quad\textrm{if } |i+k-j-l|>4 \textrm{ or } |i-k|\notin\{0,2\}.
\end{equation}
for the mass matrix. Similar results (with a similar sparsity pattern)  can be proved for stiffness matrix as well as
for the 3D case and other applications, see  \cite{Zaglmayer} for an overview.
Nevertheless, it remains to compute the nonzero entries in optimal arithmetical complexity.
The type of integrals is similar to all of the above mentioned applications. For simplicity, it is explained for $m_{ij,kl}$ in \eqref{eq3p7}
on the reference element $\refel$ with the vertices $V_1=(-1,-1)$, $V_2=(1,-1)$ and $V_3=(0,1)$. Then, one has to compute
\begin{align}\notag
m_{ij,kl}^{\refel}&=\int_{-1}^1 \int_{\frac{y-1}{2}}^{\frac{1-y}{2}} \hat{p}_{i}^{(0,0)}\left(\frac{2x}{1-y} \right)  
\hat{p}_{k}^{(0,0)}\left(\frac{2x}{1-y} \right) \left(\frac{1-y}{2} \right)^{i+k} \hat{p}^{(2i-1,0)}_j(y) \hat{p}^{(2k-1,0)}_l(y) \dx \dy \\
\overset{z=\frac{2x}{1-y}}&{=} \underbrace{\int_{-1}^{1} \hat{p}_{i}^{(0,0)}(z)  \hat{p}_{k}^{(0,0)}(z)\dz}_{:=I_{i,k}^{(0,0,0,0)} } \int_{-1}^1 \left(\frac{1-y}{2} \right)^{i+k+1} \hat{p}^{(2i-1,0)}_j(y) \hat{p}^{(2k-1,0)}_l(y) \dy .\label{DufEx}
\end{align}
The integral $I^{(0,0,0,0)}_{i,k}$ is only nonzero if $|i-k|\notin \{0,2\}$.
Inserting this into the second integral, one has to compute
\begin{equation}\label{eq3p8}
I^{(i+k+1,0,2i-1,2k-1)}_{j,l}=\int_{-1}^1 \left(\frac{1-y}{2} \right)^{i+k+1} \hat{p}^{(2i-1,0)}_j(y) \hat{p}^{(2k-1,0)}_l(y) \dy
\end{equation}
with $i,j\in \N$ and $|i-k|\leq 2$.
For the other applications, the weights in \eqref{eq3p8} differ slightly.
Our aim is to develop recursion formulas for \eqref{eq3p8} in a more general setting.

%
%
%
\section{\label{ch4}Recursion identities}
\subsection{Application: Building a bridge between special functions and FEM}
The FEM basis functions are defined on a triangle. By using the so called Duffy transformation one can transform these onto the square $(-1,1) \times (-1,1)$,
see e.g. \cite{Dubiner},\cite{Karniadakis} for details. 
Thus the problem in the FEM application boils down to one dimensional integrals over integrated Jacobi polynomials. A generalized version of one of these integrals is  by \eqref{NotIntJac} and \eqref{NotIntLeg} equivalent to 
\begin{equation*}
I_{n,m} \coloneqq \int_{-1}^1 (1-x)^\mu (1+x)^\nu P_n^{(\alpha,\beta)}(x) P_m^{(\rho,\delta)}(x) \dx,
\end{equation*}
with $n+\alpha+\beta \geq 0$ as well as $m+\rho+\delta \geq 0$.\\
We are interested in finding recursion formulas for $I_{n,m}$ with respect to $n,m$.
Our main result for the numerical application of this paper is the following theorem. 
\begin{theorem}\label{MainTheo}
 The exact value $I_{n,m} \coloneqq \int_{-1}^1 \left(\frac{1-x}{2}\right)^\mu P_n^{(\alpha,0)}(x) P_m^{(\rho,0)}(x) \dx$
 can be calculated recursively by
 \begin{align}\label{RecInt2}
  (m+n +\mu +1)~ I_{n,m} = (m+n+\alpha+\rho-\mu-1)~ I_{n-1,m-1} + (n - m  +\alpha - \mu -1)~ I_{n-1,m} + (m - n  +\rho -\mu -1)~ I_{n,m-1},
 \end{align}
 where $\alpha, \rho > -1,$ and $\mu \geq 0.$\\
Moreover the exact value $I_{n,m}^{(2)} \coloneqq \int_{-1}^1 \left(\frac{1-x}{2}\right)^\mu \phat_n^{(\alpha,0)}(x) \phat_m^{(\rho,0)}(x) \dx$
 can be calculated recursively by
 \begin{align}\label{RecInt}
  (m+n+\mu +1)~ I^{(2)}_{n,m} = (m+n+\alpha+\rho-\mu-5)~ I^{(2)}_{n-1,m-1} + (n-m+\alpha-\mu-3)~ I^{(2)}_{n-1,m} + (m - n +\rho - \mu -3)~ I^{(2)}_{n,m-1}
 \end{align}
\end{theorem}
\begin{proof}
The results are an immediate consequence of theorem \ref{Main} and corollary \ref{folg4p13}.
\end{proof}

\begin{figure}[htbp]
\begin{center}
\begin{tabular}{cc}
  \includegraphics[width=0.28\textwidth]{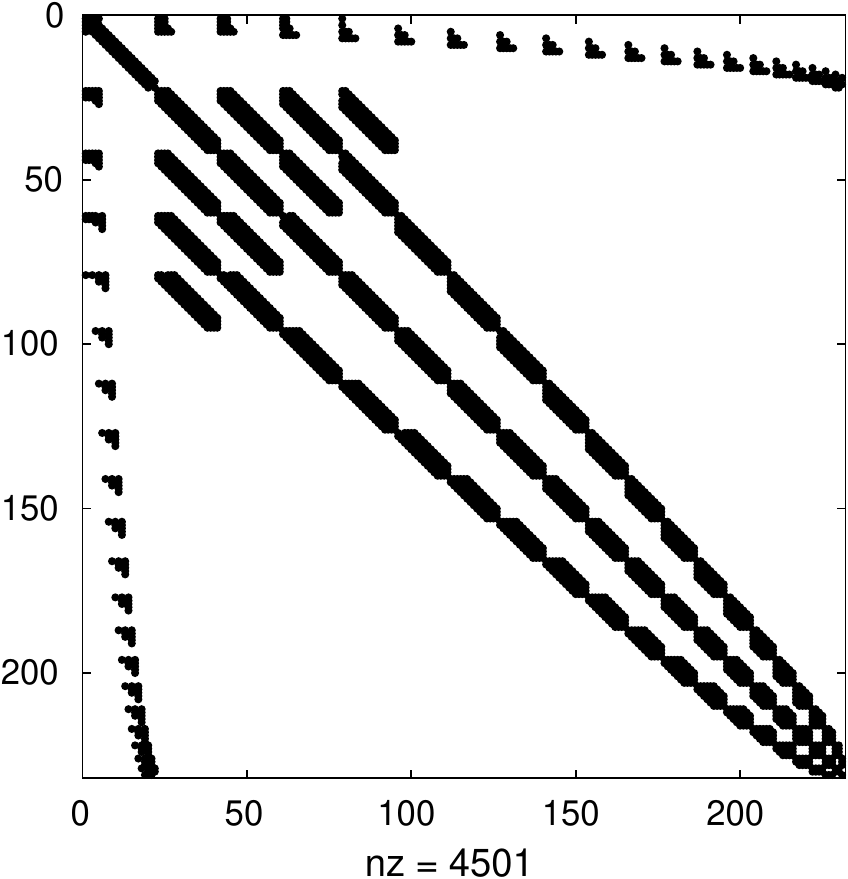}
 &
 \begin{tikzpicture}{scale=2.0}
\node [fill=red,circle, draw] (a1) at (0,6) [label=above :{$(j-1,\ell-1)$}]{};
\node [fill=red,circle, draw] (b1) at (2,6) [label=above: {$(j,\ell-1)$}]{};
\node [fill=green,circle, draw] (f1) at (2,4) [label=below :{$(j,\ell)$}]{};
\node [fill=red,circle, draw] (c1) at (0,4) [label= below :{$(j-1,\ell)$}]{};
\draw[->] (a1) -- (f1);
\draw[->] (b1) -- (f1);
\draw[->] (c1) -- (f1);
\end{tikzpicture}
\end{tabular}
\caption{\label{fig:1}Nonzero pattern of the two-dimensional mass matrix on a triangle (left), schematic representation of the recursions \eqref{RecInt}, \eqref{RecInt2} for the $(j,\ell)$ entry. }
\end{center}
\end{figure}
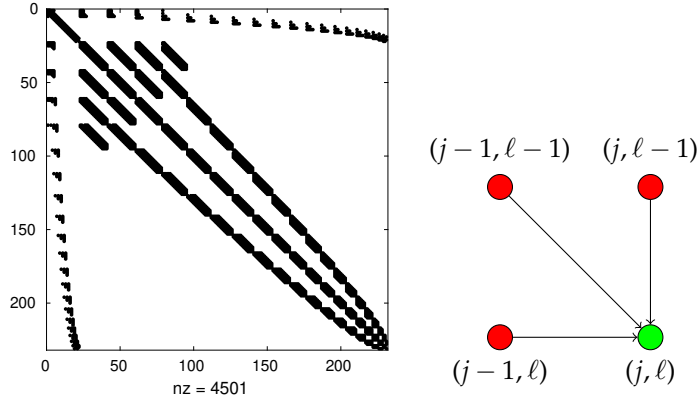

Both  recursion formulas shown in theorem \ref{MainTheo} are just special cases of the more general theorem \ref{Main}, which will be derived and proven analytically in chapter 4.4.
The right picture in figure \ref{fig:1} 
shows a schematic representation of the recursion formulas, whereas the left picture displays the typical nonzero pattern of the element matrix.
We will give now a short example, how these recursion formulas can be applied to the finite element context. 
\begin{example}
The recursion formula for the entry $I^{(2)}_{n,m}$ can be applied directly to case of the mass matrix for the $H^1$ basis functions.
On the triangle $\hat{T}$ with vertices $(-1,-1),(1,-1)$ and $(0,1)$ the interior basis functions \eqref{eq3p6} take the form  \[\Psi_{i,j} = \phat_i^{(0,0)}\left(\frac{2x}{1-y}\right) \left(\frac{1-y}{2}\right)^i \phat^{(2i-1,0)}_j (y).\]
This results in the following two one dimensional integrals, as seen above in \eqref{DufEx},
\begin{eqnarray*}
m_{ij,kl}^{\refel}\overset{z=\frac{2x}{1-y}}{=}& \underbrace{ \int_{-1}^{1} \hat{p}_{i}^{(0,0)}(z)  \hat{p}_{k}^{(0,0)}(z)\dz}_{:=I_{i,k}^{(0,0,0,0)} } \underbrace{\int_{-1}^1 \left(\frac{1-y}{2} \right)^{i+k+1} \hat{p}^{(2i-1,0)}_j(y) \hat{p}^{(2k-1,0)}_l(y) \dy.}_{:=I_{j,l}^{(i+k+1,0,2i-1,2k-1)} }
\end{eqnarray*}
Both integrals $I^{(0,0,0,0)}_{i,k}$ and $I^{(i+k+1,0,2i-1,2j-1)}_{j,l}$ can be computed using \eqref{RecInt}, though in the case of $I^{(0,0,0,0)}_{i,k}$ it reduces to
\begin{align*}
 (m+n +1)~ I^{(0,0,0,0)}_{n,m} = (m+n-5)~ I^{(0,0,0,0)}_{n-1,m-1},
\end{align*}
due to the sparsity pattern. \\
Since each recursion formula needs starting values, we compute the lowest order integrals in each block of the mass matrix, compare the left picture in figure \ref{fig:1}, by a low order quadrature. An extension to the edge functions is straight forward, since those are just special
cases of the interior functions.
\end{example}

The remainder of this paper is dedicated to proof the generalized version of \ref{MainTheo}. We start by rewriting the exact integral value.\\
By using a series representation of the Jacobi polynomial, one can write the exact value of the integral as
\begin{equation}
 I_{n,m} = \frac{ 2^{\mu+\nu+1} (\alpha+1)_n (\rho+1)_m }{n!m!} \sum_{l=0}^n \sum_{r=0}^{m}  \frac{\frac{(-n)_l (n+\alpha+\beta+1)_l}{(\alpha+1)_l l!} (-m)_r(m+\rho+\delta+1)_r}{(\rho+1)_r r!  }  \int_{-1}^1 (1-x)^{\mu+r+l} (1+x)^{\nu} \dx.
\end{equation}
The integral in the double sum is exactly Euler's Beta integral, i.e.
\begin{equation}
 \int_{-1}^1 (1-x)^{\mu+r+l} (1+x)^{\nu} \dx = B(\mu+r+l+1,\nu+1),
\end{equation}
where $B(x,y) = \frac{\Gamma(x)\Gamma(y)}{\Gamma(x+y)}$.
Thus $I_{n,m}$ can be written as
\begin{equation}\label{IForm}
 I_{n,m} = \frac{ 2^{\mu+\nu+1}  (\alpha+1)_n (\rho+1)_m}{ n! m!}  \sum_{l=0}^n  \sum_{r=0}^{m} Z_{r,l}
\end{equation}
with
\begin{equation*}
Z_{r,l} \coloneqq \frac{(-n)_l (n+\alpha+\beta+1)_l}{(\alpha+1)_l l!} \frac{(-m)_r(m+\rho+\delta+1)_r}{ (\rho+1)_r r!}  B(\mu+r+l+1,\nu+1). 
\end{equation*}
The Beta function can be expressed as using Pochhammer symbols as well, i.e.
\begin{equation*}
 B(\mu+r+l+1,\nu+1) = \frac{\Gamma(\mu+r+l+1)\Gamma(\nu+1)}{\Gamma(\mu+\nu+r+l+2)} = \Gamma(\nu+1) \frac{\Gamma(\mu+1)}{\Gamma(\mu+\nu+2)} \frac{(\mu+1)_{r+l}}{(\mu+\nu+2)_{r+l}} = B(\mu+1,\nu+1) \frac{(\mu+1)_{r+l}}{(\mu+\nu+2)_{r+l}}.
\end{equation*}

\subsection{Recurrence relation for the special case $\nu,\beta,\delta = 0$}
If we set $\nu=\beta=\delta = 0$ one can rewrite \eqref{IForm} in a generalized hypergeometric series. First, rewrite the Beta function using Pochhammer symbols, then split $r+l$ by using combinatorial arguments and use the Pfaff-Saalsch\"utz theorem \eqref{Pfaff} to reduce the double sum in \eqref{IForm}. 
After some more combinatorial arguments 
\begin{equation}
\begin{aligned}\label{Identity}
 I_{n,m} &= \int_{-1}^{1} (1-x)^\mu P^{(\alpha,0)}_n P^{(\rho,0)}_m \dx \\&=\text{const.}\, \pFq{4}{3}{-m,\rho+m+1,1,1+\mu-\alpha}{\rho +1, n+2, 1+\mu-\alpha  -n}{1} \quad \text{with } \mu \leq \alpha \text{ or } \alpha \in \R\setminus \N,
\end{aligned}
\end{equation}
follows, see \cite{Kalla} or \cite{Erdelyi} for more details.
Starting from this representation, we now prove a recursion originally obtained through the symbolic package \texttt{Guess} by Manuel Kauers \cite{Kauers}.
\begin{theorem}\label{EasyRec}
 Let $\alpha = \rho$, $\mu= \beta =\delta =\nu = 0$. Then the recursion relation
 \begin{equation}\label{recursion}
  c_0 I_{n+1,m+1} + c_1 I_{n+1,m} + c_2 I_{n,m+1} + c_3 I_{n,m} = 0,
 \end{equation}
 with the coefficients
 \begin{align*}
 c_0 &= (m+n+3),\\
 c_1 &= -(\alpha+m-n-1),\\
 c_2 &= -(\alpha-m+n-1),\\ 
 c_3 &= -(2\alpha + m +n +1)
 \end{align*}
 holds.
\end{theorem}
\begin{proof}
   We start by deriving a contiguous relation of the hypergeometric series in \eqref{Identity}.
  Since 
  \begin{equation*}
  \pFq{4}{3}{-m,\alpha+m+1,1,1-\alpha}{\alpha +1, n+2, 1-\alpha  -n}{1} = \sum_{\kappa=0}^{\infty} \underbrace{\frac{(-m)_\kappa (m+\alpha+1)_\kappa (1)_\kappa (1-\alpha)_\kappa}{(\alpha + 1)_\kappa (n+2)_\kappa (1-\alpha -n)_\kappa} \frac{1}{\kappa !}}_{\eqqcolon \phi_\kappa}
  \end{equation*}
  the contiguous functions can be written as
  \begin{dmath*}
  \pFq{4}{3}{-m-1,\alpha+m+2,1,1-\alpha}{\alpha +1, n+2, 1-\alpha  -n}{1} = \sum_{\kappa = 0}^\infty \frac{(-m-1)(\alpha+m+1+\kappa)}{(\kappa-m-1)(\alpha+m+1)} \phi_\kappa
  \end{dmath*}
  and
  \begin{dmath*}
  \pFq{4}{3}{-m,\alpha+m+1,1,1-\alpha}{\alpha +1, n+3, -\alpha  -n}{1} = \sum_{\kappa = 0}^\infty \frac{(2+n)(\kappa-\alpha-n)}{(2+n+\kappa)(-n-\alpha)} \phi_\kappa.
  \end{dmath*}
  Hence, in order to find a recursion of the form \eqref{recursion},
  \begin{equation*}
   0 = \sum_{\kappa = 0}^{\infty} \left(- c_0 \frac{(m+\alpha+1+\kappa)(\kappa-\alpha-n)}{(2+n+\kappa)(m-1+\kappa)} - c_1 \frac{(m+\alpha+1+\kappa)}{(-m-1+\kappa)} - c_2 \frac{(\kappa - \alpha -n)}{(2+n+\kappa)} + c_3 \right) \phi_\kappa
  \end{equation*}
needs to hold. Expanding the fractions in the same denominator reduces the condition to
\begin{dmath*}
 -c_0 (m+1+\alpha+\kappa)(\kappa - \alpha - n) - c_1 (m+1+a+\kappa)(2+n+\kappa) - c_2 (\kappa-\alpha-n)(\kappa-m-1) + c_1 (\kappa-m-1)(n+2+\kappa) = 0.
\end{dmath*}
To find the coefficients $c_i, i= 0,\dots,3$, we view this as polynomial of $\kappa$ and equate the coefficients of $\kappa$ to zero, leading to the equations
\begin{align*}
&-c_0 -c_1 -c_2 +c_3 &= 0,\\
& -c_0 (k+1-n) - c_1 (-\alpha-n-m-1) - c_2 (m+n+\alpha+3) + c_3(n+1-m)&= 0,\\
& -c_0(m+1+\alpha)(-\alpha-n) - c_1(-\alpha-n)(-m-1) -c_2(m+1+\alpha)(n+2) + c_3(-m-1)(n+2) &= 0.
\end{align*}
This leads to the underdetermined system
\begin{equation*}
 \begin{pmatrix}
  -(m+1+\alpha)(-\alpha-n) & -(-\alpha-n)(-m-1) & -(m+1+\alpha)(n+2) & (-m-1)(n+2)\\-(m+1-n) & -(-\alpha-n-m-1) & -(m+n+\alpha+3) & (n+1-m)\\ -1&-1&-1&1 
 \end{pmatrix}
\begin{pmatrix}
 c_0\\c_1\\c_2\\c_3
\end{pmatrix}
= \bf{0},
\end{equation*}
then we can choose $c_0 = 1$ and bring this to the right-hand side
\begin{equation*}
  \begin{pmatrix}
   -(-\alpha-n)(-m-1) & -(m+1+\alpha)(n+2) & (-m-1)(n+2)\\ -(-\alpha-n-m-1) & -(m+n+\alpha+3) & (n+1-m)\\ -1&-1&1
 \end{pmatrix}
\begin{pmatrix}
 c_1\\c_2\\c_3
\end{pmatrix}
= \begin{pmatrix}
   -(m+1+\alpha)(-\alpha-n) \\-(m+1-n)\\ -1
  \end{pmatrix}
\end{equation*}
Solving this system leads to
\begin{align*}
 c_1 = \frac{-(\alpha+m-n-1)}{(m+n+3)},\\
 c_2 = \frac{-(\alpha-m+n-1)}{(m+n+3)},\\ 
 c_3 = \frac{-(2\alpha + m +n +1)}{(m+n+3)}.
 \end{align*}
 Rescaling $c_0$ leads to the proposed recurrence relation.
\end{proof}\\
Further simplification of the generalized hypergeometric series in \eqref{Identity} is not that that trivial. There is a summation formula for a balanced ${}_4F_3$ in (\cite{Rakha1}, \cite{Rakha2}), but this can not be applied here since the coefficients of the series don't match.
A summation by using Whipple's transformation and Dougall's summation (see e.g. \cite{Askey}) works only for the case $n=m$, otherwise the resulting transformed series isn't well-poised.\\
Moreover the more general case $\mu, \nu, \beta$ and $\delta$ arbitrarily can't be represented by a ${}_4F_3$ since neither of the sums in \eqref{IForm} is summable by the Pfaff-Saalsch\"utz theorem. To be precise, 
those general series are $\nu+1$ balanced. We therefore use a more general approach.

\subsection{Kamp\'{e} de F\'{e}riet Series}
The concept of hypergeometric series can be extended to the multivariate case. Examples are the Appell series \cite{Appell}, the Kamp\'{e} de F\'{e}riet series\cite{Appell} or the Lauricella series.
\begin{definition}
For $p_1,p_2,p_3,q_1,q_2,q_3 \in \N$ and coefficients $(a_1, \dots a_{q_1}),(b_1, \dots b_{p_1}), \dots$ arbitrarily the series
 \begin{equation*}
  F^{p_1;p_2;p_3}_{q_1;q_2;q_3} = \sum_{n,m = 0}^\infty \frac{ \prod_{i=1}^{p_1} (k_i)_{n+m} \prod_{i=1}^{p_2} (a_i)_n \prod_{i=1}^{p_3} (b_i)_m}{ \prod_{i=1}^{q_1} (l_i)_{n+m} \prod_{i=1}^{q_2} (c_i)_n \prod_{i=1}^{q_3} (d_i)_m} \frac{x^n y^m}{n! m!}
 \end{equation*}
is called Kamp\'{e} de F\'{e}riet series.\\
\end{definition}
The notation is due to Burchnall and Chaundy (see \cite{Burchnall},\cite{Burchnall2}). More information, in particular on the convergence theory of such series, can be found e.g. in \cite{Exton}.\\
We can write $I_{n,m}$ now as a generalized Kamp\'{e} de F\'{e}riet series, i.e.
\begin{equation}
\begin{aligned}\label{KFSeriesEasy}
 F \coloneqq I_{n,m} = 2^{\mu+\nu+1} \frac{(\alpha+1)_n (\rho+1)_m \Beta(\nu+1,\mu+1) }{n!m!} F^{1;2;2}_{1;1;1}\left(\begin{matrix}\mu+1 & ; &-n&n+\alpha+\beta+1&;&-m&m+\rho+\delta+1 \\ \mu+\nu+2&;&\alpha+1& &;& \rho+1&\end{matrix} ;1;1 \right),
\end{aligned}
\end{equation}
where $\Beta(x,y) = \frac{\Gamma(x)\Gamma(y)}{\Gamma(x+y)}$ is the usual Beta-function as above. $F$ converges, since it is terminating due to the coefficients $-n$ and $-m$.\\
Furthermore we omit indices of $F$ to keep the notation as simple as possible.
We denote the forward shift of parameters by one analogously to earlier, i.e.,
\begin{equation*}
\begin{aligned}
F(\alpha+) &= 2^{\mu+\nu+1} \frac{(\alpha+2)_n (\rho+1)_m \Beta(\nu+1,\mu+1) }{ n! m!} \, F^{1;2;2}_{1;1;1}\left(\begin{matrix}\mu+1 & ; &-n&n+\alpha+\beta+2&;&-m&m+\rho+\delta+1 \\ \mu+\nu+2&;&\alpha+2& &;& \rho+1&\end{matrix} ;1;1 \right)\\
F(n+) &= 2^{\mu+\nu+1} \frac{ (\alpha+1)_n (\rho+1)_m \Beta(\nu+1,\mu+1)}{ (n+1)!m!}\, F^{1;2;2}_{1;1;1}\left(\begin{matrix}\mu+1 & ; &-n-1&n+\alpha+\beta+2&;&-m&m+\rho+\delta+1 \\ \mu+\nu+2&;&\alpha+1& &;& \rho+1&\end{matrix} ;1;1 \right)\\
 F(\mu+) &= 2^{\mu+\nu+2} \frac{(\alpha+1)_n (\rho+1)_m \Beta(\nu+1,\mu+2) }{n! m!}\, F^{1;2;2}_{1;1;1}\left(\begin{matrix}\mu+2 & ; &-n&n+\alpha+\beta+1&;&-m&m+\rho+\delta+1 \\ \mu+\nu+2&;&\alpha+1& &;& \rho+1&\end{matrix} ;1;1 \right)\\
 F(\nu+)& = 2^{\mu+\nu+2} \frac{(\alpha+1)_n (\rho+1)_m \Beta(\nu+2,\mu+1) }{n! m!}\, F^{1;2;2}_{1;1;1}\left(\begin{matrix}\mu+1 & ; &-n&n+\alpha+\beta+1&;&-m&m+\rho+\delta+1 \\ \mu+\nu+2&;&\alpha+1& &;& \rho+1&\end{matrix} ;1;1 \right)\\
 &\cdots 
\end{aligned}
 \end{equation*}
\subsection{The general case}\label{RecurrenceRelation}
We can now directly derive recurrence relations for \eqref{KFSeriesEasy}
by using the recurrence relations of Lemma \ref{le:Rek} above. Namely

\begin{corollary}
 \begin{dmath*}
  (\alpha+\beta+n) F = (\beta + n) F(\beta-) + (\alpha+n) F(\alpha-)
 \end{dmath*}

 \begin{dmath}\label{Rec2}
  -\frac{1}{2}(2+\alpha+\beta+2n) F(\alpha+,\mu+) = (n+1) F(n+) - (1+\alpha+n) F
 \end{dmath}

 \begin{dmath}\label{Rec3}
  \frac{1}{2}(2+\alpha+\beta+2n) F(\beta+,\nu+) = (n+1) F(n+) + (1+\beta+n) F
 \end{dmath}

 \begin{dmath}\label{Rec4}
  (\alpha+\beta+2n) F(\beta-) = (\alpha+\beta+n) F + (\alpha+n) F(n-)
 \end{dmath}

 \begin{dmath}\label{Rec5}
  (\alpha+\beta+2n) F(\alpha-) = (\alpha+\beta+n) F - (\beta+n) F(n-)
 \end{dmath}

 \begin{dmath*}
  2 F = F(\nu+,\beta+) + F(\mu+,\alpha+)
 \end{dmath*}

 \begin{dmath*}
  F(n-) = F(\beta-) - F(\alpha-)
 \end{dmath*}
\end{corollary}
and 7 analogue recurrence formulas in $m,\rho$ and $\delta$, where
 \begin{dmath}\label{Rec4_2}
  (\rho+\delta+2m) F(\delta-) = (\rho+\delta+m) F + (\rho+m) F(m-)
 \end{dmath}
 is one of them.\\
Furthermore we can derive more recurrence relations by linear combinations of the above.
\begin{corollary}
 The following recurrence formula holds
 \begin{align}\label{basicrec}
  -2(1+\alpha+n) F - (1+\beta+n)F(\alpha+) + (2+\alpha+\beta+2n)F(\alpha+,\mu+) - (\alpha+\beta+1) F(n+) +(2+\alpha+\beta+n) F(n+,\alpha+) =0
 \end{align}
\end{corollary}
\begin{proof}
 The equation is just a linear combination of \eqref{Rec2} and \eqref{Rec5}.
Start by transforming \eqref{Rec2} and \eqref{Rec5}
\begin{align*}
 0 &= (2+\alpha+\beta+2n) F(\alpha+, \mu+) - 2(n+1) F(n+) - 2(1+\alpha+n) F,\\
 0 &= -(\alpha+\beta+2n+3) F(n+) - (\beta+n+1)F(\alpha+) + (\alpha +\beta +n +2) F(n+,\alpha+).
\end{align*}
Adding these equations lead to~\eqref{basicrec}.
 \end{proof}\\
One important, but rather trivial relation is given by
\begin{equation}\label{munu}
 2F = F(\nu+) + F(\mu+),
\end{equation}
which follows from $2 P_n^{(\alpha,\beta)}= (1+x)P_n^{(\alpha,\beta)} +(1-x)P_n^{(\alpha,\beta)}$.

The following relations will later be proven by using a generalized form of \eqref{KFSeriesEasy} in the appendix. Set $x = y = 1$ in lemma \ref{lem1} and \ref{lem2}, then the following corollary holds  
\begin{corollary}
 \begin{align}
  \label{MixedRec1}(n+m+\mu+\nu+4) F(n+,m+,\nu+) &= (\alpha+n+1)F(m+,\beta+,\nu+) + (\rho + m +1) F(n+,\delta+,\nu+) + 2 (\nu+1) F(n+,m+)\\
  \label{MixedRec2}(n+\alpha+\beta+m+\rho+\delta-\mu-\nu+1)F &= (n+\alpha+\beta+1) F(\beta+) + (m+\rho+\delta+1)F(\delta+) + 2 \nu F(\nu-)
 \end{align}
\end{corollary}
Since the weights in the Jacobi polynomials are interchangeable, see \eqref{atob}, the following two relations can be derived as well.
\begin{corollary}
 \begin{align}
  \label{MixedRec3}(n+m+\mu+\nu+4) F(n+,m+,\mu+) &= -(\beta+n+1)F(m+,\alpha+,\mu+) - (\delta + m +1) F(n+,\rho+,\mu+) + 2 (\mu+1) F(n+,m+)\\
  \label{MixedRec4}(n+\alpha+\beta+m+\rho+\delta-\mu-\nu+1)F &= (n+\alpha+\beta+1) F(\alpha+) + (m+\rho+\delta+1)F(\rho+) + 2 \mu F(\mu-),
 \end{align}
\end{corollary}
see lemma \ref{lem3} and \ref{lem4} and set $x=y=1$.
\paragraph{5 - point recurrence relation}
There are some known starlike recurrence relations, see \cite{PillweinPamm}. Those can be derived in this context as follows. 
\begin{corollary}
 We have two mixed recurrence relations
 \begin{dmath}\label{Mixed1}
  (2m+\rho+\delta+1) \left( (n+1) F(n+,\alpha-) - (\alpha+n) F(\alpha-) \right) = (2n+ \alpha +\beta +1) \left( (m+1)F(m+,\rho-) - (\rho+m)F(\rho-)\right)
 \end{dmath}
 \begin{dmath}\label{Mixed2}
  (2m+\rho+\delta+1) \left( (n+1) F(n+,\beta-) - (\beta+n) F(\beta-) \right) = (2n + \alpha +\beta +1) \left( (m+1)F(m+,\delta-) - (\delta+m)F(\delta-)\right)
 \end{dmath}
\end{corollary}
\begin{proof}
 Take \eqref{Rec2} and replace $\alpha$ by $\alpha-1$ to derive
 \begin{dmath*}
  F(\mu+) = \frac{-2}{2n+\alpha+\beta+1} \left((n+1) F(n+,\alpha-) - (\alpha+n)F(\alpha-)\right)
 \end{dmath*}
and $\rho$ by $\rho-1$ to derive
\begin{dmath*}
 F(\mu+) = \frac{-2}{2m+\rho+\delta+1} \left((m+1) F(m+,\rho-) - (\rho+m)F(\rho-)\right).
 \end{dmath*}	
Setting both right hand sides equal yields \eqref{Mixed1}. The second mixed relation follows analogously from the recursion formula~\eqref{Rec3}.
\end{proof}\\
The mixed relations \eqref{Mixed1} and \eqref{Mixed2} yield some 5-point recurrence relations with support $(m,n),(m-1,n),(m+1,n),(m,n-1),(m,n+1),$  see also \cite{PillweinPamm}. 
\begin{theorem}
\begin{equation}\label{Mixed5Point}
 \begin{aligned} (2m+\rho+\delta)_3 \left[ (n+1)\left((n+\alpha+\beta+1)(2n+\alpha+\beta) F(n+) + (n+\alpha+1) F\right) + (\beta+n)\left( (n+\alpha+\beta) F + (n + \alpha) F(n-)\right)\right]\\ 
 = (2n+\alpha+\beta)_3 \left[ (m+1)\left((m+\rho+\delta+1)(2m+\rho+\delta) F(m+) + (m+\rho+1) F\right) + (m+\delta)\left( (m+\rho+\delta) F + (m + \delta) F(m-)\right)\right]
 \end{aligned}
 \end{equation}
\begin{equation}\label{Mixed5Point2}
 \begin{aligned} (2m+\rho+\delta)_3 \left[ (n+1)\left((n+\alpha+\beta+1)(2n+\alpha+\beta) F(n+) - (n+\beta+1) F\right) - (\alpha+n)\left( (n+\alpha+\beta) F - (n + \beta) F(n-)\right)\right]\\ 
 = (2n+\alpha+\beta)_3 \left[ (m+1)\left((m+\rho+\delta+1)(2m+\rho+\delta) F(m+) - (m+\delta+1) F\right) + (m+\rho)\left( (m+\rho+\delta) F + (m + \rho) F(m-)\right)\right]
 \end{aligned}\end{equation}
  \begin{align*} (2m+\rho+\delta)_3 \left[ (n+1)\left(2(n+\alpha+\beta+1)(2n+\alpha+\beta) F(n+) + (\alpha-\beta) F\right) + \left( (\beta-\alpha)(n+\alpha+\beta) F + (n+\alpha)(n + \beta) F(n-)\right)\right]\\ 
 = (2n+\alpha+\beta)_3 \left[ (m+1)\left(2(m+\rho+\delta+1)(2m+\rho+\delta) F(m+) + (\rho-\delta) F\right) + \left( (\delta-\rho)(m+\rho+\delta) F + (m + \rho)(m+\delta) F(m-)\right)\right]
 \end{align*}
\end{theorem}
\begin{proof}
 Take the first mixed relation \eqref{Mixed1} and replace all series by \eqref{Rec5}, this yields the first equation. The second equation follows by using \eqref{Mixed2} with \eqref{Rec4}. Lastly, the third equation can be derived by linear combination of \eqref{Mixed5Point} and \eqref{Mixed5Point2}
\end{proof} \\
Alternatively one can use the 3-term recursion \eqref{3term} to proof the same result as in \cite{PillweinPamm}.
\paragraph{Recurrence relation}
Multiple recurrence relations similar to \eqref{recursion} can be proven. 
\begin{lemma}
 Let $F = I_{n,m}$, where $I_{n,m}$ is as in \eqref{KFSeriesEasy}. Then the following recurrence relation holds
 \begin{dmath}\label{recnu}
  (n+\alpha+\beta+1)(m+\rho+\delta+1)\left((n+m+\mu+\nu+4) F(n+,m+,\nu+)-2(\nu+1)F(n+,m+)\right) = (\alpha+n+1)(m+\rho+\delta+1)\left( (n+\alpha+\beta-m-\mu-\nu-2)F(m+,\nu+) + 2(\nu+1)F(m+)\right) + (\rho+m+1)(n+\alpha+\beta+1)\left((-n+m+\rho+\delta-\mu-\nu-2)F(n+,\nu+) + 2(\nu+1)F(n+)\right)+(\rho+m+1)(\alpha+n+1)\left((n+\alpha+\beta+m+\rho+\delta-\mu-\nu)F(\nu+) + 2(\nu+1) F\right). 
 \end{dmath}
\end{lemma}
\begin{proof}
 Start with equation \eqref{MixedRec1}
 \begin{equation*}
  (n+m+\mu+\nu+4) F(n+,m+,\nu+) - 2 (\nu+1) F(n+,m+)= (\alpha+n+1)F(m+,\beta+,\nu+) + (\rho + m +1) F(n+,\delta+,\nu+). 
 \end{equation*}
Replace both terms of the RHS by using shifted versions of equation \eqref{MixedRec2}, i.e.
\begin{equation}\begin{aligned}\label{FirstReplace}
 &\frac{(n+\alpha+1)}{(n+\alpha+\beta+1)} (n+\alpha+\beta+1) F(m+,\beta+,\nu+) \\ &=\frac{(n+\alpha+1)}{(n+\alpha+\beta+1)}\left[(n+\alpha+\beta+m+\rho+\delta -\mu -\nu +1)F(m+,\nu+) - (m+\rho+\delta+2) F(m+,\delta+,\nu+) + 2(\nu+1) F(m+)\right],
\end{aligned}\end{equation}
and 
\begin{equation}\begin{aligned}\label{SecondReplace}
 &\frac{(m+\rho+1)}{(m+\rho+\delta+1)} (m+\rho+\delta+1) F(n+,\delta+,\nu+) \\ &=\frac{(m+\rho+1)}{(m+\rho+\delta+1)}\left[(n+\alpha+\beta+m+\rho+\delta -\mu -\nu +1)F(n+,\nu+) - (n+\alpha+\beta+2) F(n+,\beta+,\nu+) + 2(\nu+1) F(n+)\right].
\end{aligned}\end{equation}
Moreover use the shifted relation \eqref{Rec4} for the middle part of the last two equations, i.e.
\begin{align}
 \label{ThirdReplace}(m+\rho+\delta+2) F(m+,\delta+,\nu+) = (2m+\rho+\delta+3) F(m+,\nu+) - (m+\rho+1) F(\delta+,\nu+)\\
 \label{FourthReplace}(n+\alpha+\beta+2) F(n+,\beta+,\nu+) = (2n+\alpha+\beta+3) F(n+,\nu+) - (n+\alpha+1) F(\beta+,\nu+).
\end{align}
Lastly replace the remaining terms by a shifted version of \eqref{MixedRec2},
\begin{align*}
 &\frac{(n+\alpha+1)}{(n+\alpha+\beta+1)}(m+\rho+1) F(\delta+,\nu+) + \frac{(m+\rho+1)}{(m+\rho+\delta+1)}(n+\alpha+1) F(\beta+,\nu+)\\
 &= \frac{(n+\alpha+1)(m+\rho+1)}{(n+\alpha+\beta+1)(m+\rho+\delta+1)} \left( (n+\alpha+\beta+1) F(\beta+,\nu+) + (m+\rho+\delta+1)F(\delta+,\nu+)\right)\\
 &= \frac{(n+\alpha+1)(m+\rho+1)}{(n+\alpha+\beta+1)(m+\rho+\delta+1)} \left( (n+\alpha+\beta+m+\rho+\delta-\mu-\nu) F(\nu+) + 2(\nu+1) F\right)
\end{align*}
The claim follows from combining the above with the remaining terms of \eqref{FirstReplace},  \eqref{SecondReplace}, \eqref{ThirdReplace} and \eqref{FourthReplace}. 
\end{proof}\\

Since $\alpha$ and $\beta$ or $\rho$ and $\delta$ are interchangeable, the following lemma can be proven by using \eqref{MixedRec3} and \eqref{MixedRec4} instead of \eqref{MixedRec1} and \eqref{MixedRec2}. Hence
\begin{lemma}
  \begin{dmath}\label{recmu}
  (n+\alpha+\beta+1)(m+\rho+\delta+1)\left((n+m+\mu+\nu+4) F(n+,m+,\mu+)-2(\mu+1)F(n+,m+)\right) = -(\beta+n+1)(m+\rho+\delta+1)\left( (n+\alpha+\beta-m-\mu-\nu-2)F(m+,\mu+) + 2(\mu+1)F(m+)\right) - (\delta+m+1)(n+\alpha+\beta+1)\left((-n+m+\rho+\delta-\mu-\nu-2)F(n+,\mu+) + 2(\mu+1)F(n+)\right)+(\delta+m+1)(\beta+n+1)\left((n+\alpha+\beta+m+\rho+\delta-\mu-\nu)F(\mu+) + 2(\mu+1) F\right). 
 \end{dmath}
\end{lemma}

Both of these recursion formulas have the drawback, that terms with $\nu+1$ or $\mu+1$ vanish only for $\nu = -1$ or $\mu = -1,$ which correspond to the special cases $(1+x)^0$ or $(1-x)^0.$ If the steps of the proof are slightly adjusted, a recursion formula, which applies to more cases, can be proven. The following theorem is the main result of this paper:
\begin{theorem}\label{Main}
 Let $F = I_{n,m}$, where $I_{n,m}$ is as in \eqref{KFSeriesEasy}. Then the following recurrence relation holds
  \begin{dmath}\label{HardRec}
  (n+1)(m+1)\left((n+m+\mu+\nu+4) F(n+,m+,\nu+)-2(\nu+1-\beta-\delta)F(n+,m+)\right) = (n+\beta+1)(m+1)\left( (n+\alpha+\beta-m-\mu-\nu-2)F(m+,\nu+) + 2(\nu+1-\beta-\delta)F(m+)\right) + (n+1)(m+\delta+1)\left((-n+m+\rho+\delta-\mu-\nu-2)F(n+,\nu+) + 2(\nu+1-\beta-\delta)F(n+)\right)+(n+\beta+1)(m+\delta+1)\left((n+\alpha+\beta+m+\rho+\delta-\mu-\nu)F(\nu+) + 2(\nu+1-\beta-\delta) F\right). 
 \end{dmath}
\end{theorem}
\begin{proof}
 Again start with recursion \eqref{MixedRec1}, i.e.
  \begin{equation*}
  (n+m+\mu+\nu+4) F(n+,m+,\nu+) - 2 (\nu+1) F(n+,m+)= (\alpha+n+1)F(m+,\beta+,\nu+) + (\rho + m +1) F(n+,\delta+,\nu+), 
 \end{equation*}
 now add $2(\beta+\delta)F(n+,m+)$ to both sides and multiply by the factor $(n+1)(m+1)$ on both sides. Thus
   \begin{equation*}
   \begin{aligned}
  (n+1)(m+1)&\left((n+m+\mu+\nu+4) F(n+,m+,\nu+) -2 (\nu+1-\beta-\delta) F(n+,m+)\right) \\ &= (n+1)(m+1)\left((\alpha+n+1)F(m+,\beta+,\nu+) + (\rho + m +1) F(n+,\delta+,\nu+) + 2(\beta+\delta)F(n+,m+)\right) \\ &=\operatorname{RHS}. 
 \end{aligned}
 \end{equation*}
 Instead of multiplying with $1$, as in the proof for \eqref{recnu}, we will add a $0$ to expand the RHS.
 Hence
 \begin{equation*}
 \begin{aligned}
  \operatorname{RHS}=&(n+1)(m+1)\Big((n+\alpha+\beta+1)F(m+,\beta+,\nu+) + (m+\rho +\delta+1) F(n+,\delta+,\nu+) \\&-(\beta F(m+,\beta+,\nu+)+\delta F(n+,\delta+,\nu+))+2(\beta+\delta)F(n+,m+)\Big)\\
  =&(n+\beta+1)(m+1)(n+\alpha+\beta+1)F(m+,\beta+,\nu+) + (n+1)(m+\delta+1)(m+\rho+\delta+1)F(n+,\delta+,\nu+)\\ &- \beta (m+1)(n+\alpha+\beta+1) F(m+,\beta+,\nu+) - \delta (n+1)(m+\rho+\delta+1) F(n+,\delta+,\nu+)\\
  &- \beta (n+1)(m+1) F(m+,\beta+,\nu+) - \delta (n+1)(m+1) F(n+,\delta+,\nu+) + 2(n+1)(m+1)(\beta+\delta) F(n+,m+).
 \end{aligned}
 \end{equation*}
 After adding up the additional terms recurrence relation \eqref{Rec3} can be used. This gives
\begin{equation*}
\begin{aligned}
\operatorname{RHS} = &(n+\beta+1)(m+1)(n+\alpha+\beta+1)F(m+,\beta+,\nu+) + (n+1)(m+\delta+1)(m+\rho+\delta+1)F(n+,\delta+,\nu+)\\ &-2\beta(m+1)\left((n+1)F(n+,m+)+(n+\beta+1) F(m+)\right) - 2\delta(n+1)\left((m+1)F(n+,m+) + (m+\delta+1) F(n+)\right) \\&+ 2(n+1)(m+1)(\beta+\delta) F(n+,m+1)\\
 = &(n+\beta+1)(m+1)(n+\alpha+\beta+1)F(m+,\beta+,\nu+) + (n+1)(m+\delta+1)(m+\rho+\delta+1)F(n+,\delta,\nu+) \\&-2\beta(m+1)(n+\beta+1) F(m+) - 2\delta(n+1)(m+\delta+1) F(n+).
\end{aligned}
 \end{equation*}
Now use the mixed relation \eqref{MixedRec2} 
\begin{equation*}
 \begin{aligned}
   \operatorname{RHS} = &(n+\beta+1)(m+1)\left[(n+\alpha+\beta+m+\rho+\delta-\mu-\nu+1)F(m+,\nu+) +2(\nu+1) F(m+) + (m+\rho+1)F(\delta+,\nu+)\right]\\ &+(n+1)(m+\delta+1)\left[(n+\alpha+\beta+m+\rho+\delta-\mu-\nu+1)F(n+,\nu+) +2(\nu+1) F(n+) + (n+\alpha+1)F(\beta+,\nu+)\right] \\&-2\beta(m+1)(n+\beta+1) F(m+) - 2\delta(n+1)(m+\delta+1) F(n+)\\
   = &(n+\beta+1)(m+1)\left[(n+\alpha+\beta+m+\rho+\delta-\mu-\nu+1)F(m+,\nu+) +2(\nu+1-\beta-\delta) F(m+) + (m+\rho+1)F(\delta+,\nu+)\right]\\ &+(n+1)(m+\delta+1)\left[(n+\alpha+\beta+m+\rho+\delta-\mu-\nu+1)F(n+,\nu+) +2(\nu+1-\beta-\delta) F(n+) + (n+\alpha+1)F(\beta+,\nu+)\right]\\ &+ 2 \delta (n+\beta+1)(m+1)F(m+) + 2\beta (n+1)(m+\delta+1)F(n+).
 \end{aligned}
\end{equation*}
Consider only a part of the RHS to shorten the notation. Begin by transforming $F(n+)$ or $F(m+)$ back to the form $F(\beta+,\nu+)$ or $F(\delta+,\nu+)$ by equation \eqref{Rec3}, i.e.
\begin{equation*}
\begin{aligned}
(n+\beta+1)&(m+1)(m+\rho+1)F(\delta+,\nu+) + (n+1)(m+\delta+1)(n+\alpha+1)F(\beta+,\nu+)\\ &+ 2\delta(n+\beta+1)(m+1)F(m+) + 2\beta(n+1)(m+\delta+1)F(n+)\\
=&(n+\beta+1)(m+1)(m+\rho+1) F(\delta+,\nu+) + (n+1)(m+\delta+1)(n+\alpha+1)F(\beta+,\nu+)\\ &+ \delta (n+\beta+1)(2m+\rho+\delta+2) F(\delta+,\nu+) + \beta (m+\delta+1)(2n +\alpha+\beta+2) F(\beta+,\nu+)\\ &+ 2\delta (n+\beta+1)(m+\delta+1) F + 2 \beta (n+\beta+1)(m+\delta+1) F\\
=& (m+1)(n+\beta+1)(m+\rho+\delta+1) F(\delta+,\nu+) + (n+1)(m+\delta+1)(n+\alpha+\beta+1)F(\beta+,\nu+)\\ &+ \delta(n+\beta+1)(m+\rho+\delta+1) F(\delta+,\nu+) + \beta(m+\delta+1)(n+\alpha+\beta+1)F(\beta+,\nu+)\\&+ 2 \beta (n+\beta+1)(m+\delta+1) F\\
=& (m+\delta+1)(n+\beta+1)(m+\rho+\delta+1) F(\delta+,\nu+) + (n+\beta+1)(m+\delta+1)(n+\alpha+\beta+1)F(\beta+,\nu+)\\&+ 2 \beta (n+\beta+1)(m+\delta+1) F.
\end{aligned}
\end{equation*}
In the last step use again \eqref{MixedRec2}, then the claim follows.
\end{proof}\\
Setting $\nu = -1$, $\beta = \delta = 0$ and $\alpha = \rho$ in \eqref{HardRec} annihilates the coefficient $(\nu + 1 - \beta - \delta)$ and thus reduces \eqref{HardRec} to
\begin{equation*}
\begin{aligned}
  (n+1)(m+1)\left((n+m+3) F(n+,m+)\right) =& (n+1)(m+1)\left( (n+\alpha-m-1)F(m+)\right)\\ &+ (n+1)(m+1)\left((-n+m+\alpha-1)F(n+)\right)+(n+1)(m+1)\left((n+m+2\alpha+1)F\right), 
\end{aligned}
\end{equation*}
which is again theorem \ref{EasyRec}.
\begin{corollary}\label{folg4p13}
If $\nu + 1 = \beta +\delta$ equation \eqref{HardRec} reduces to
\begin{dmath}
  (n+1)(m+1)\left((n+m+\mu+\nu+4) F(n+,m+,\nu+)\right) = (n+\beta+1)(m+1)\left( (n+\alpha+\beta-m-\mu-\nu-2)F(m+,\nu+)\right) + (n+1)(m+\delta+1)\left((-n+m+\rho+\delta-\mu-\nu-2)F(n+,\nu+)\right)+(n+\beta+1)(m+\delta+1)\left((n+\alpha+m+\rho-\mu-1)F(\nu+)\right). 
 \end{dmath}
\end{corollary}
\begin{remark}
 The case $\nu+1 = \beta + \delta$ is corresponds to the integrated Jacobi polynomials $\widehat{p}_n^{(\alpha,0)}(x) = \frac{1-x}{n} P_{n-1}^{(\alpha-1,1)}(x),$ where $\beta = \delta = 1, \nu = 1, n \geq 2, \alpha > 1.$
\end{remark}
\section{\label{ch5}Numerical aspects}
Usually an arbitrary high order finite element matrix is computed by Gaussian quadrature. The steps of an element assembly algorithm can be summarized in 3 steps:
\begin{enumerate}
 \item Calculate the $1$D quadrature points and weights: $\mathcal{O}(n)$, see e.g. \cite{TownsendGauss}.
 \item Evaluate Jacobi polynomials or more generally the basis functions: $\mathcal{O} (n)$ in $1$D for $n$ quadrature points, see e.g. \cite{Higham}.
 \item Gaussian quadrature and in higher dimension sum factorization:$\mathcal{O}(n^{d+1})$, see \cite{Orszag}, \cite{Gerdes}.
\end{enumerate}
The ansatz by using the recursion formulas, skips the first two steps and replaces the third by a summation, which is independent of $n$. For the ease of representation we provide a $1$D example.
\subsection{1D - Example}
\begin{figure}
\begin{subfigure}[b]{0.48 \textwidth}
  \includegraphics[width = 0.68\textwidth]{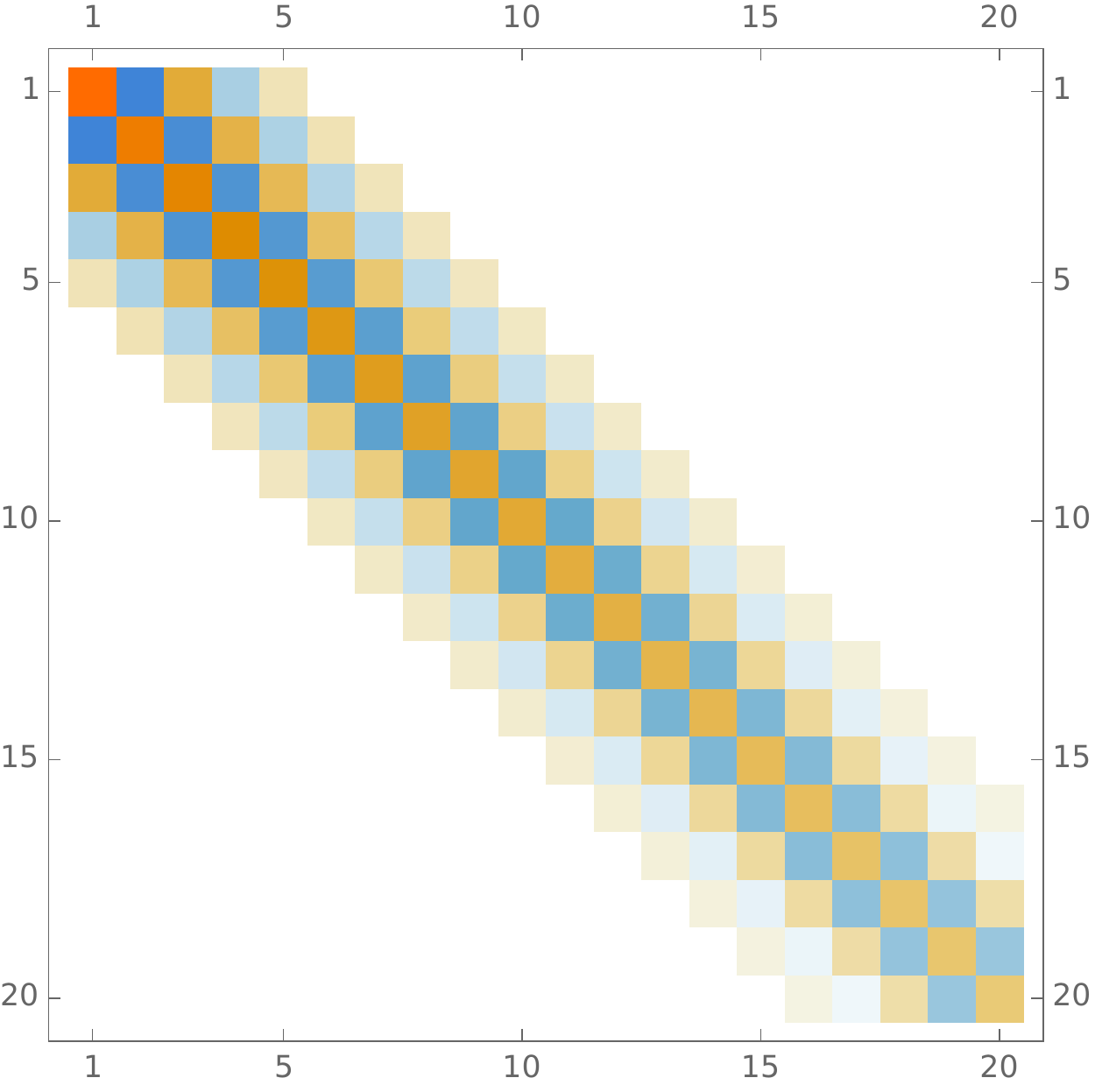}
  \caption{\label{SparsityP}Sparsity pattern of the example}
 \end{subfigure}
\begin{subfigure}[b]{0.48 \textwidth}
 \includegraphics[width = \textwidth]{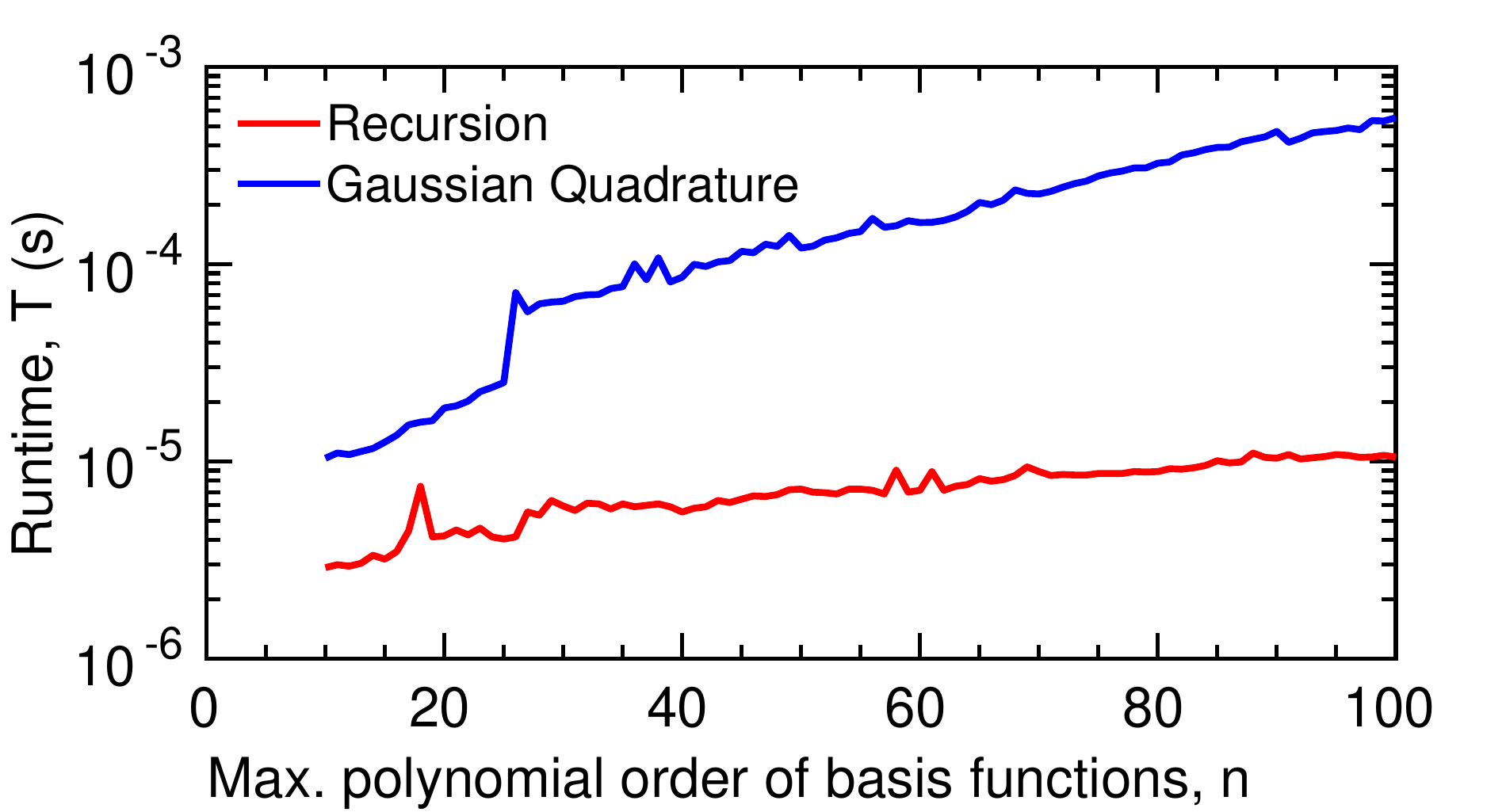}
 \caption{\label{Runtime}Comparing the runtime of Gaussian quadrature with recursion formulas}
\end{subfigure}
\caption{\label{Gram}}
\end{figure}

Let $\scalarp{f(x)}{g(x)} = \int_{-1}^{1} f(x) g(x) \left(\frac{(1-x)}{2}\right)^8 \dx$ be a weighted scalar product. Let $\phi_{n}(x) = P_{n}^{(4,0)}(x)$, then calculate the Gram matrix $G = \left[G_{i,j}\right]_{i,j}^{n_{max}^2} =\left[\scalarp{\phi_i}{\phi_j}\right]_{i,j}^{n_{max}^2}.$
The resulting sparsity pattern can be seen in figure \ref{Gram}\subref{SparsityP}. To compare the standard assembly routine with the recursive version, we assume now, that the quadrature points, weights and the basis functions are tabulated, 
i.e. we only need to perform step 3 of the standard assembly routine.
In figure \ref{Gram}\subref{Runtime} the runtime of a Matlab implementation of both assembly routines are compared. We measured the assembly time of the Gram matrix for the total polynomial order $10 < n_{max} < 160$ for the integration of each non zero value $\scalarp{\phi_i}{\phi_j}.$
As expected the recursive version is a lot faster than the Gaussian quadrature, even in $1$D.  While this recurrence relations hold even for low order polynomials, the main strength is the application to higher order or spectral methods, where
two integrals in $2D$ or three in $3D$ are computed for each non zero entry. The number of flops per entry remains constant in the recursive case, thus we achieved optimal arithmetical complexity for each integral dimension.

\begin{remark}
Polynomial coefficients $q_n(x)$ of the degree $n$ in the integrand, i.e. $\scalarp{q_n(x) \phi_i(x)}{\phi_j(x)}$ can be handled as well. Since the Legendre polynomials are a basis of $\mathcal{P}_n,$ we can write $q_n(x) = \sum_{i=0}^n c_i L_i$.
The following two corollaries show, that we are able to connect integrals of three orthogonal polynomials and hence we are able to compute polynomial coefficients for the integrand as well, though depending on $q_n$ a lot more flops are needed.
\begin{corollary}
With given $a\in\N$ the integrand 
\[
 I_{a,n,m} = \left(\frac{1-x}{2}\right)^{i+j+1} L_a(x) P_n^{(2i,0)}(x) P_m^{(2j,0)}(x)
\]
satisfies the following recursion formula
\begin{align*}
 (2 + a + i + j + m + n) I_{a,n,m} &= -(2 - a - i - j - m - n) I_{n-1, m-1, a-1} - (2 + a - i - j - m - n) I_{n-1, m-1,a}  \\
  &+ (-2 + a + i - j - m + n) I_{n-1, m,a-1} - (2 + a - i + j + m - n) I_{n-1, m,a}\\ &+ (-2 + a - i + j + m - n) I_{n, m-1,a-1}
  - (2 + a + i - j - m + n) I_{n, m-1, a}\\ &+ (-2 + a - i - j - m - n) I_{n, m, a-1} 
\end{align*}

\end{corollary}
A similar result holds for the integrated Jacobi polynomials.
\begin{corollary}
 For the integrand 
\[
 I_{a,n,m} = \left(\frac{1-x}{2}\right)^{i+j+1} L_a(x) \phat_n^{2i,0}(x) \phat_m^{2j,0}(x)
\]
holds the following recursion formula
\begin{align*}
 (2 + a + i + j + m + n) I_{a,n,m} &= -(6 - a - i - j - m - n) I_{n-1, m-1, a-1} - (6 + a - i - j - m - n) I_{n-1, m-1,a}  \\
  &+ (-4 + a + i - j - m + n) I_{n-1, m,a-1} - (4 + a - i + j + m - n) I_{n-1, m,a}\\&+ (-4 + a - i + j + m - n) I_{n, m-1,a-1}
  - (4 + a + i - j - m + n) I_{n, m-1, a}\\ &+ (-2 + a - i - j - m - n) I_{n, m, a-1} 
\end{align*}

\end{corollary}
\end{remark}
\section{Conclusion}
We were able to derive and prove generalized versions of recurrence formulas, which were originally derived by symbolic computation.
One of the main complications in high order finite elements is the assembly of local mass and stiffness matrices. Using these newly computed recurrences one is able to gain a 
better asymptotic complexity than the state of the art method, sum factorization. Furthermore this approach does not only skip the numerical quadrature, but the initialization step as well. 
From the special functions point of view the representation as Kamp\'{e} de F\'{e}riet series in itself is interesting. 

\section{Acknowledgement}
The first author has been supported by the Deutsche Forschungsgemeinschaft (DFG) under Germany’s
Excellence Strategy within the Cluster of Excellence PhoenixD (EXC 2122, Project ID 390833453).
\newpage
\begin{appendix}
\section{Derivation of mixed recurrence relation}
\subsection{Generalization}
More general one can look at the Kamp\'{e} de F\'{e}riet series where $x$ and $y$ are not equal to 1, i.e.
\begin{equation}\label{KFSeries}
 F := F^{1;2;2}_{1;1;1}\left(\begin{matrix}\mu+1 & ; &-n&n+\alpha+\beta+1&;&-m&m+\rho+\delta+1 \\ \mu+\nu+2&;&\alpha+1& &;& \rho+1&\end{matrix} ;x;y \right). 
\end{equation}
Furthermore we omit the indices of $F$ to keep the notation as simple as possible. 
\subsubsection{Recurrence formula}
One can extend the proofs for the contiguous relations of a hypergeometric series to the general case of the Kamp\'{e} de F\'{e}riet series using the techniques in Rainville \cite{Rainville} Ch.4.
Let $Z$ be a generalization of $F$ to arbitrary coefficients, i.e.
\begin{equation*}
Z = \sum_{n=0}^\infty \sum_{m=0}^\infty \underbrace{\frac{(a)_n (b)_n (f)_m (g)_m (d)_{n+m} x^n y^m}{n! m! (c)_n (h)_m (e)_{n+m}} }_{\eqqcolon \tau_{n} \tau_{m} \tau_{n+m}}
\end{equation*}
be a general Kamp\'{e} de F\'{e}riet series. It has similar contiguous functions to the six contiguous functions of the Gaussian hypergeometric series, namely
\begin{align*}
Z(a+) = \sum_{n=0}^\infty \frac{a+n}{a} \tau_{n}\tau_{m}\tau_{n+m}, \quad Z(a-) = \sum_{n=0}^\infty \frac{a-1}{a-1+n} \tau_{n}\tau_{m}\tau_{n+m},\\
Z(b+) = \sum_{n=0}^\infty \frac{b+n}{b} \tau_{n}\tau_{m}\tau_{n+m}, \quad Z(b-) = \sum_{n=0}^\infty \frac{b-1}{b-1+n} \tau_{n}\tau_{m}\tau_{n+m},\\
Z(c+) = \sum_{n=0}^\infty \frac{c}{c+n} \tau_{n}\tau_{m}\tau_{n+m}, \quad Z(c-) = \sum_{n=0}^\infty \frac{c-1+n}{c-1} \tau_{n}\tau_{m}\tau_{n+m},
 \end{align*}
obviously one can find 6 similar functions for $f,g$ and $h$. In the following omit $(f)_m,(g)_m$ and $(h)_m$ to simplify the notation.\\
Use the differential operator $\theta_x = x \frac{\partial}{\partial x}$. This leads to
\begin{equation*}
 (\theta_x +a) Z = \sum_{n=0}^\infty (a+n)\tau_{n}\tau_{m}\tau_{n+m}
\end{equation*}
and hence
\begin{equation}
\begin{aligned}\label{CR}
 (\theta_x +a) Z &= a Z(a+),\\
 (\theta_x +c -1) Z &= (c-1) Z(c-),
 \end{aligned}
\end{equation}
and so on. Following one of the calculations in \cite{Rainville}, one can derive
\begin{equation*}
 \begin{aligned}
  \theta_x Z(a-) &= \sum_{n=1}^\infty \sum_{m=0}^\infty \frac{(a-1)_n (b)_n (d)_{n+m} \tau_{m} x^n y^m}{(c)_n (n-1)! (e)_{n+m}} \\
  &= x \sum_{n=0}^\infty \sum_{m=0}^\infty \frac{(a-1)_{n+1} (b)_{n+1} (d)_{n+m+1} \tau_{m} x^n y^m}{(c)_{n+1} n! (e)_{n+m+1}}\\
  &= x (a-1) \frac{d}{e} \sum_{n=0}^\infty \sum_{m=0}^\infty \frac{(b+n)(d+n+m)}{(c+n)(e+n+m)} \tau_{n+m} \tau_n \tau_{m}\\
  &= x(a-1) \frac{d}{e} Z(d+,e+) - \frac{(a-1)(c-b) d}{c e} x Z(c+,d+,e+),
 \end{aligned}
\end{equation*}
since $ \frac{b+n}{c+n} = 1- \frac{c-b}{c+n} \frac{c}{c}$. Now replace $\theta_x Z(a-1)$ by \eqref{CR} with $a$ replaced by $(a-1)$.
Thus 
\begin{equation*}
 \frac{e}{d} Z = \frac{e}{d} Z(a-) + x Z(d+,e+) - \frac{(c-b)}{c} x Z(c+,d+,e+),
\end{equation*}
since $a$ and $b$ are interchangeable, there holds 
\begin{equation*}
 \frac{e}{d} Z = \frac{e}{d} Z(b-) + x Z(d+,e+) - \frac{(c-a)}{c} x Z(c+,d+,e+).
\end{equation*}
Subtracting the last two equations from each other yields
\begin{equation*}
 0 = e (Z(a-)-Z(b-) +c^{-1} d (b-a) x Z(c+,d+,e+)
\end{equation*}
and if we set $b$ to $b+1$
\begin{equation*}
 0 = e (Z(a-,b+)-Z) +c^{-1} d (b+1-a) x Z(b+,c+,d+,e+).
\end{equation*}
Setting $a=-n, b= n+\alpha+\beta+1, c= \alpha+1$ and $x=1$ (and the respective values for $f,g,h$) leads after simplification to recursion formula \eqref{Rec2}.

\subsubsection{Proof of mixed relations}
Moreover the recurrence relations, which we have seen in section \ref{RecurrenceRelation},  can be derived for the more general case \eqref{KFSeries}. 
All of the following recursion hold for $x=1$ and $y=1$, which follows just from the recurrence formula of the Jacobi polynomials, as we have seen.\\  
To prove the mixed relations \eqref{MixedRec1} and \eqref{MixedRec2} consider 
\begin{equation}\begin{aligned}\label{KFSeries2}
F \coloneqq &\frac{ 2^{\mu+\nu+1}  (\alpha+1)_n (\rho+1)_m \Beta(\nu+1,\mu+1)}{ n! m!}\\ &F^{1;2;2}_{1;1;1}\left(\begin{matrix}\mu+1 & ; &-n&n+\alpha+\beta+1&;&-m&m+\rho+\delta+1 \\ \mu+\nu+2&;&\alpha+1& &;& \rho+1&\end{matrix} ;x;y \right). 
\end{aligned}\end{equation}
Denote the contiguous functions as usual, i.e.
\begin{equation*}
\begin{aligned}
F(\alpha+) &= \frac{ 2^{\mu+\nu+1}  (\alpha+2)_n (\rho+1)_m \Beta(\nu+1,\mu+1)}{ n! m!}F\left(\begin{matrix}\mu+1 & ; &-n&n+\alpha+\beta+2&;&-m&m+\rho+\delta+1 \\ \mu+\nu+2&;&\alpha+2& &;& \rho+1&\end{matrix} ;x;y \right)\\
F(n+) &= \frac{ 2^{\mu+\nu+1}  (\alpha+1)_{n+1} (\rho+1)_m \Beta(\nu+1,\mu+1)}{ (n+1)! m!}F\left(\begin{matrix}\mu+1 & ; &-n-1&n+\alpha+\beta+2&;&-m&m+\rho+\delta+1 \\ \mu+\nu+2&;&\alpha+1& &;& \rho+1&\end{matrix} ;x;y \right)\\
 F(\mu+) &= \frac{ 2^{\mu+\nu+2}  (\alpha+1)_n (\rho+1)_m \Beta(\nu+1,\mu+2)}{ n! m!}F\left(\begin{matrix}\mu+2 & ; &-n&n+\alpha+\beta+1&;&-m&m+\rho+\delta+1 \\ \mu+\nu+2&;&\alpha+1& &;& \rho+1&\end{matrix} ;x;y \right)\\
 F(\nu+)& = \frac{ 2^{\mu+\nu+2}  (\alpha+1)_n (\rho+1)_m \Beta(\nu+2,\mu+1)}{ n! m!}F\left(\begin{matrix}\mu+1 & ; &-n&n+\alpha+\beta+1&;&-m&m+\rho+\delta+1 \\ \mu+\nu+2&;&\alpha+1& &;& \rho+1&\end{matrix} ;x;y \right)\\
 &\cdots 
\end{aligned}
 \end{equation*}
 \begin{lemma}
Let $\theta_x = x\frac{\partial}{\partial x}$ and $\theta_y = y \frac{\partial}{\partial y}$, then the following differential equations hold
\begin{align}
 \label{dx1}(\theta_x - n) F &= -(n+\alpha) F(n-,\beta+)\\
 \label{dx2}(\theta_x + n +\alpha+\beta+1) F &= (n+\alpha+\beta+1) F(\beta+)\\
 \label{dy1}(\theta_y - m) F &= -(m+\rho) F(m-,\delta+)\\
 \label{dy2}(\theta_y +m +\rho +\delta +1) F &= (m+\rho+\delta+1) F(\delta+)\\
 \label{dxdy}(\theta_x +\theta_y +\mu+\nu+1) F &= 2\nu F(\nu-)
\end{align}
\end{lemma}
\begin{proof}
For $F$ as in \eqref{KFSeries2}
\begin{equation*}\begin{aligned}
(\theta_x - n) F = &\frac{ 2^{\mu+\nu+1}  (\alpha+1)_n(\rho+1)_m \Beta(\nu+1,\mu+1)}{ n! m!}\\ &\sum_{k = 0}^\infty \sum_{l = 0}^\infty \frac{ (\mu+1)_{k+l} (-n)_k (n+\alpha+\beta+1)_k (-m)_l (m+\rho+\delta+1)_l }{(\mu+\nu+2)_{k+l}(\alpha+1)_k (\rho+1)_l} (\theta_x -n) \frac{x^k y^l}{k! l!},\\
 =&\frac{ 2^{\mu+\nu+1}  (\alpha+1)_n (\rho+1)_m \Beta(\nu+1,\mu+1)}{ n!m!}\\ &\sum_{k = 0}^\infty \sum_{l = 0}^\infty \frac{ (\mu+1)_{k+l} (-n)_k (n+\alpha+\beta+1)_k (-m)_l (m+\rho+\delta+1)_l }{(\mu+\nu+2)_{k+l}(\alpha+1)_k (\rho+1)_l} (k-n) \frac{x^{k-1} y^l}{k! l!}.
 \end{aligned}\end{equation*}
 As usual replace  $(k-n) (-n)_k$ by $(-n) (-n+1)_k$, thus
 \begin{equation*}\begin{aligned}
  =&\frac{ (-n) 2^{\mu+\nu+1}  (\alpha+1)_n(\rho+1)_m \Beta(\nu+1,\mu+1)}{ n!m!}\\ &\sum_{k = 0}^\infty \sum_{l = 0}^\infty \frac{ (\mu+1)_{k+l} (-n+1)_k (n+\alpha+\beta+1)_k (-m)_l (m+\rho+\delta+1)_l }{(\mu+\nu+2)_{k+l} (\alpha+1)_k (\rho+1)_l} \frac{x^k y^l}{k! l!}\\
  =&\frac{ 2^{\mu+\nu+1} (\alpha+1)_n(\rho+1)_m \Beta(\nu+1,\mu+1) }{ n!m!}\\ &F^{1;2;2}_{1;1;1}\left(\begin{matrix}\mu+1 & ; &-n-1& (n-1) +\alpha+\beta+2&;&-m&m+\rho+\delta+1 \\ \mu+\nu+2&;&\alpha+1& &;& \rho+1&\end{matrix} ;x;y \right),
\end{aligned}\end{equation*}
since we changed $n$ to $n-1$ in the series, we need to equate $n+\alpha+\beta+1$ as well. This is done by raising $\beta$ by one. Furthermore we need to change the parts in prefactor, to accommodate $n-1$ as well. Since
\begin{equation}
 \begin{aligned}
  \frac{ -n 2^{\mu+\nu+1}  (\alpha+1)_n(\rho+1)_m \Beta(\nu+1,\mu+1)}{ n!m!} = - \frac{ 2^{\mu+\nu+1} (\rho+1)_m \Beta(\nu+1,\mu+1) (\alpha +n) (\alpha+1)_{n-1}}{m! (n-1)!},
 \end{aligned}
\end{equation}
the equation \eqref{dx1} follows. Relation \eqref{dy1} follows analogue to it. Relation \eqref{dx2} or \eqref{dy2} follow directly by applying the differential operators, since the prefactor doesn't need to be changed.\\
For the last relation \eqref{dxdy} the differential operator $(\theta_x +\theta_y +\mu+\nu+1)$ applied to $x^k y^l$ yields $(k+l+\mu+\nu+1)$. This reduces $(\nu+\mu+2)_{k+l}$ in the denominator to $(\nu+\mu+2)_{k+l-1}$. Therefore we multiply the series with $\frac{\nu+\mu+1}{\nu+\mu+1},$ such that the part in the denominator becomes $(\nu+\mu+1)_{k+l}$.
Since it is only a change in the denominator it is rather a change in $\nu$ than in $\mu$. The rest follows using a property of the Beta-function, i.e.
\begin{equation*}
 B(x+1,y) = B(x,y) \frac{x}{x+y}.
\end{equation*}
\end{proof}\\

Subtracting \eqref{dx1} and \eqref{dy1} from \eqref{dxdy} proves \eqref{MixedRec1}.
\begin{lemma}\label{lem1}
 For $F$ as in \eqref{KFSeries2} the following contiguous recurrence relation hold,
 \begin{equation*}
  (n+m+\mu+\nu+4) F(n+,m+,\nu+) = (n+\alpha+1) F(m+,\beta+,\nu+) + (m+\rho+1) F(n+,\delta+,\nu+) + 2(\nu+1) F(\nu+)
 \end{equation*}
\end{lemma}
 Similar \eqref{MixedRec2} can be proven. Take \eqref{dx2} and \eqref{dy2} and subtract \eqref{dxdy}
 \begin{lemma}\label{lem2}
  For $F$ as in \eqref{KFSeries2} the following contiguous recurrence relation hold,
 \begin{equation*}
  (n+\alpha+\beta+m+\rho+\delta-\mu-\nu+1) F = (n+\alpha+\beta+1) F(\beta+) + (m+\rho+\delta+1) F(\delta+) + 2 \nu F(\nu-)
 \end{equation*}
 \end{lemma}
 The easiest way of deriving \eqref{MixedRec3} and \eqref{MixedRec4} is to introduce another formulation for $F$. Recall that Jacobi polynomials can be expressed as
 \begin{equation*}
  P^{(\alpha,\beta)}_n = \frac{(-1)^n (\beta+1)_n}{n!} \pFq{2}{1}{-n,n+\alpha+\beta+1}{\beta+1}{\frac{1+x}{2}}
 \end{equation*}
as well. Using this expression in the derivation of the Kamp\'{e} de F\'{e}riet series yield the analogue form
\begin{equation}\begin{aligned}\label{KFSeriesAna}
 \widetilde{F} \coloneqq &\frac{ (-1)^{n+m} 2^{\mu+\nu+1} (\beta+1)_n (\delta+1)_m \Beta(\nu+1,\mu+1) }{m! n!} \\&F^{1;2;2}_{1;1;1}\left(\begin{matrix}\nu+1 & ; &-n&n+\alpha+\beta+1&;&-m&m+\rho+\delta+1 \\ \mu+\nu+2&;&\beta+1& &;& \delta+1&\end{matrix} ;x;y \right).
\end{aligned}\end{equation}
Then the differential equations can be derived on the same way as seen before. Thus
\begin{lemma}
\begin{align}
 \label{dx3}(\theta_x - n) \widetilde{F} &= (n+\beta) \widetilde{F}(n-,\alpha+),\\
 \label{dx4}(\theta_x + n +\alpha+\beta+1) \widetilde{F} &= (n+\alpha+\beta+1) \widetilde{F}(\alpha+),\\
 \label{dy3}(\theta_y - m) \widetilde{F} &= (m+\delta) \widetilde{F}(m-,\rho+),\\
 \label{dy4}(\theta_y +m +\rho +\delta +1) \widetilde{F} &= (m+\rho+\delta+1) \widetilde{F}(\rho+),\\
 \label{dxdy2}(\theta_x +\theta_y +\mu+\nu+1) \widetilde{F} &= 2\mu \widetilde{F}(\mu-).
\end{align}
\end{lemma}
Now \eqref{MixedRec3} follows by subtracting \eqref{dx3} and \eqref{dy3} from \eqref{dxdy2}
\begin{lemma}\label{lem3}
 For $F$ as in \eqref{KFSeries2} or \eqref{KFSeriesAna} the following contiguous recurrence relation hold,
 \begin{equation*}
  (n+m+\mu+\nu+4) F(n+,m+,\mu+) = -(n+\beta+1) F(m+,\alpha+,\mu+) - (m+\delta+1) F(n+,\rho+,\mu+) + 2(\mu+1) F(\mu+)
 \end{equation*}
\end{lemma}
and analogue
\begin{lemma}\label{lem4}
 For $F$ as in \eqref{KFSeries2} or \eqref{KFSeriesAna} the following contiguous recurrence relation hold,
 \begin{equation*}
  (n+\alpha+\beta+m+\rho+\delta-\mu-\nu+1) F = (n+\alpha+\beta+1) F(\alpha+) + (m+\rho+\delta+1) F(\rho+) + 2 \mu F(\mu-).
 \end{equation*}
\end{lemma}

\begin{remark}
 Following Burchnall and Chaundy \cite{Burchnall}, \cite{Burchnall2} one can easily compute an expansion base, i.e. $F = \sum_{i=0}^\infty A_i(x) A_i(y),$ for $F.$ 
\end{remark}

\end{appendix}
 \bibliographystyle{siam}
\bibliography{literatur,buch}
\end{document}